\documentclass[12pt]{article}

\usepackage{hyperref}
\usepackage{amssymb}
\usepackage{amsmath}
\usepackage{graphics}
\usepackage{epsfig}
\newtheorem{claim}{\bf \t}[part]

\numberwithin{Assumption}{section} \numberwithin{Corollary}{section}
\numberwithin{Definition}{section} \numberwithin{equation}{section}
\numberwithin{Example}{section} \numberwithin{Lemma}{section}
\numberwithin{Proposition}{section} \numberwithin{Remark}{section}
\numberwithin{Theorem}{section}

\def\v{\varepsilon}

\def\x{\xi}
\def\t{\theta}
\def\T{\Theta}
\def\k{\kappa}

\def\m{\mu}

\def\g{\gamma}\def\G{\Gamma}
\def\d{\delta}
\def\l{\lambda}
\def\f{\frac}
\def\p{\phi}
\def\r{\rho}

\def\s{\sigma}

\def\di{\displaystyle}
\def\i{\infty}

\def\text#1{{\rm #1}}
\begin{document}

\title{\Large \bf Large-time Behavior of Solutions to the Inflow Problem of Full Compressible
 Navier-Stokes  Equations }
\author{\small \textbf{Xiaohong Qin}\thanks{Department of Mathematics, Nanjing
University of Science and Technology, Nanjing 210094, China. X. Qin
is supported in part by NSFC (grant No. 10901083).
 E-mail: xqin@amss.ac.cn.}
\qquad\qquad  \textbf{Yi Wang}\thanks{Institute of Applied
Mathematics, AMSS, CAS, Beijing 100190, China. Y. Wang is supported
by NSFC (grant No. 10801128). E-mail: wangyi@amss.ac.cn. }}
\date{}
\maketitle
\begin{abstract}
 Large-time behavior of solutions to the inflow problem  of full
  compressible Navier-Stokes equations is investigated on the half line $\mathbf{R}_+=(0,+\infty)$. The wave structure which contains
four waves: the transonic(or degenerate) boundary
  layer solution, 1-rarefaction wave, viscous  2-contact wave and
  3-rarefaction wave to the inflow problem is described  and the asymptotic stability of the
  superposition of the above four wave patterns to the inflow problem  of full
  compressible Navier-Stokes equations is proven under some smallness conditions. The proof is given by the elementary energy analysis
   based on the underlying wave structure. The main points
  in the proof are the degeneracies of the transonic boundary layer solution and the wave
  interactions in the superposition wave.
\end{abstract}


 {\bf AMS SC2000:} 35L60, 35L65

\section{Introduction }
\setcounter{equation}{0} \setcounter{Assumption}{0}
\setcounter{Theorem}{0} \setcounter{Proposition}{0}
\setcounter{Corollary}{0} \setcounter{Lemma}{0}

In this paper, we consider an initial-boundary-value problem for
full compressible Navier-Stokes equations in  \emph{Eulerian}
coordinates on the half line $\mathbf{R}_+=(0,+\infty)$
\begin{eqnarray}\label{(1.1)}
\begin{cases}
    \rho   _{t}+ (\rho u) _{x} = 0,   \cr
    (\rho u)_{t}+  \big {(} \rho u^2 +p\big {)}_x
    =(\mu  u_{x})_x, &x>0,~t>0,\cr
    \left[  \rho\left( e+\frac{1}{2}u^2\right)\right ] _t
    +\left[   \rho u\left( e+\frac{1}{2}u^2\right)+ pu \right ]_x
    = (\kappa\theta_{ x}+\mu uu_x )_x
\end{cases}
\end{eqnarray}
where $\rho (  t, x )> 0 $, $u(t, x ) $, $\theta (t, x )> 0 $, $p(t,
x )> 0 $ and $e(t, x )> 0  $   represent the mass density, the
velocity, the absolute temperature, the pressure, and the specific
internal energy of the gas respectively and $\mu
>0$ is the coefficient of viscosity, $\kappa>0$ is the coefficient of  heat
conduction. Here we assume that both $\mu$ and $\kappa$ are positive
constants. Let $v=\frac{1}{\rho}(>0)$ and $s$ denote the specific
volume and the entropy of the gas, respectively. Then by the second
law of thermodynamics, we have for the ideal polytropic gas
\begin{eqnarray}
p= R  v^{ -1 }  \theta  =A v^{ -\gamma }\exp \left(\frac{\g-1}{ R}s
\right) ,~~~e(v, \theta)=\frac{R}{\gamma-1} \theta,~~\label{(1.2)}
\end{eqnarray}
where $\gamma >1$ denotes the adiabatic exponent of gas, and $A$ and
$R$ are positive constants.

We consider the initial-boundary-value problem \eqref{(1.1)} with
the initial values
\begin{equation}
  (\r,u, \theta)( 0, x)=(\r_0,u_0,\t_0)(x)\rightarrow (\r_{+},u_{+},\theta_{+}) ~~ \text{as} ~~ x \rightarrow
 +\infty,~~
\inf\limits_{x\in \mathbf{R}_+} (\r_0,  \theta_0)(x)
>0\label{(1.3)}
\end{equation}
where $\r_{+}>0$,  $  u_{+} $  and $\theta_{+}>0 $ are given
constants.

As pointed out by \cite{Matsumura}, the boundary conditions to the
half space problem (\ref{(1.1)}) can be proposed as one of the
following three cases:

 \noindent Case I. outflow problem (negative
velocity on the boundary):
$$
 u(t, x)|_{x=0}=u_-<0, ~~\t(t, x)|_{x=0}=\t_-. \eqno(1.4)_1
$$
\noindent Case II. impermeable wall problem (zero velocity on the
boundary):
$$
u(t, x)|_{x=0}=0, ~~\t(t, x)|_{x=0}=\t_-.\eqno(1.4)_2
$$
\noindent Case III. inflow problem (positive velocity on the
boundary):
$$
u(t, x)|_{x=0}=u_->0, ~~\r(t, x)|_{x=0}=\r_-,~~\t(t,
 x)|_{x=0}=\t_-.\eqno(1.4)_3
$$
 Here all the $\r_{-}>0$,  $  u_{-} $ and
$\theta_{-}>0 $ in (1.4) are prescribed constants and of course we
assume that the initial values \eqref{(1.3)} and the boundary
conditions (1.4) satisfy the compatibility condition at the origin.
Notice that in Cases I and II, the density $\r_-$ on the boundary
$\{x=0\}$ could not be given, but in Case III, $\r_-$ must be
imposed due to the well-posedness theory of the hyperbolic equation
$\eqref{(1.1)}_1$.

In the present paper, we are concerned with the large-time behavior
of the solutions to the inflow problem (Case III) of the full
compressible Navier-Stokes equations (1.1), (1.3) and $(1.4)_3$. The
large-time behavior of the solutions to the compressible
Navier-Stokes equations \eqref{(1.1)} is closely related to the
corresponding Euler system
\setcounter{equation}{4}
\begin{equation}
\begin{cases}
  \r_t+(\r u)_x=0,\cr   (\r u)_t+\big(\r u^2+p\big)_x=0,\cr
 \big[\r\big(e+\f{u^2}{2}\big)\big]  _t+\big[\r u\big(e+\f{u^2}{2}\big)+pu\big]_x=0.
\end{cases}
 \label{(1.5)}
\end{equation}
 The Euler system (\ref{(1.5)}) is a typical example of the
hyperbolic conservation laws. It is well-known that the main feature
of the solutions to the hyperbolic conservation laws is the
formation of the shock wave no matter how smooth the initial values
are. The Euler system (\ref{(1.5)}) contains three basic wave
patterns, that is, two nonlinear waves, called shock wave and
rarefaction wave and one linear wave called contact discontinuity in
the solutions to the Riemann problem. The above three dilation
invariant wave solutions and their linear superpositions in the
increasing order of characteristic speed, i.e., Riemann solutions,
govern both local and large-time behavior of solutions to the Euler
system and so govern the large-time behavior of the solutions to the
compressible Navier-Stokes equations \eqref{(1.1)}.

There have been a large amount of literature on the large-time
behavior of solutions to the Cauchy problem of the compressible
fluid system \eqref{(1.1)} towards the viscous version of the basic
wave patterns. We refer to \cite{Duan-Liu-Zhao},
\cite{Huang-Li-Matsumura}, \cite{Huang-Matsumura-Xin},
\cite{Huang-Xin-Yang}, \cite{Kawashima-Matsumura}, \cite{Liu},
\cite{Liu-Xin-1}, \cite{Liu-Xin-2}, \cite{Matsumura-Nishihara-1},
\cite{Nishihara-Yang-Zhao}, \cite{Szepessy-Xin}, \cite{Xin} and some
references therein. All these works show that the large-time
behavior of the solutions to the
 Cauchy problem   is basically governed by the Riemann
solutions   to its corresponding hyperbolic system.

Recently, the initial-boundary value problem of (\ref{(1.1)})
attracts increasing interest  because it has more physical meanings
and of course produces some new mathematical difficulties due to the
boundary effect. Not only basic wave patterns but also a new wave,
which is called  boundary layer solution (BL-solution for brevity)
\cite{Matsumura}, may appear in the IVBP case. Matsumura
\cite{Matsumura} proposes a criterion on the question when the
BL-solution forms to the isentropic Navier-Stokes equations, where
the entropy of the gas is assumed to be constant and the equation
$\eqref{(1.1)}_3$ for the energy conservation is neglected. The
argument in \cite{Matsumura} for the isentropic Navier-Stokes
equations can also be applied to the full Navier-Stokes equations
(\ref{(1.1)}), see \cite{Huang-Li-Shi} for details. Consider the
Riemann problem to the Euler equations (\ref{(1.5)}), where the
initial right state of the Riemann data  is given by the far field
state $(\r_+,u_+,\t_+)$ in (\ref{(1.3)}), and the left end state
$(\r_-,u_-,\t_-)$ is given by the all possible states which are
consistent with the boundary condition (1.4) at $\{x=0\}$. Note that
to the outflow problem, $\r_-$ can not be prescribed and is free on
the boundary. On one hand, when the left end state is uniquely
determined so that the value at the boundary $\{x=0\}$ of the
solution to the Riemann problem is consistent with the boundary
condition, we expect that no BL-solution occurs. On the other hand,
if the value of the solution to the Riemann problem on the boundary
is not consistent with the boundary condition for any admissible
left end state, we expect a BL-solution which compensates the gap
comes up. Such BL-solution could be constructed by the stationary
solution to Navier-Stokes equations. The existence and stability of
the BL-solution (to the inflow or outflow problems, to the
isentropic or full Navier-Stokes equations) are studied extensively
by many authors, see \cite{Huang-Li-Shi},
\cite{Huang-Matsumura-Shi}, \cite{Huang-Qin},
\cite{Kawashima-Nishibata-Zhu}, \cite{Matsumura}
\cite{Matsumura-Nishihara-3}, \cite{Qin-Wang}, \cite{Zhu}, etc.

Now we review some recent works on the large-time behavior of the
solutions to the inflow problem  of the full Naiver-Stokes equation
(\ref{(1.1)}), (\ref{(1.3)}), $(1.4)_3$ by Huang-Li-Shi
\cite{Huang-Li-Shi} and Qin-Wang \cite{Qin-Wang}. In
\cite{Qin-Wang}, we rigorously prove the existence (or
non-existence) of BL-solution to the inflow problem (\ref{(1.1)}),
(\ref{(1.3)}), $(1.4)_3$ when the right end state
$(\rho_+,u_+,\theta_+)$ belongs to the subsonic, transonic and
supersonic regions respectively. When $(\r_\pm,u_\pm,\t_\pm)$ both
belong to the subsonic region, the BL-solution is expected and the
stability of this BL-solution and its superposition with the
3-rarefaction wave is proved under some smallness assumptions in
\cite{Huang-Li-Shi}. The stability of the superposition of the
subsonic BL-solution, the viscous 2-contact wave and 3-rarefaction
wave is shown in \cite{Qin-Wang} under the condition that the
amplitude of BL-solution and the contact wave is small enough but
the amplitude of the rarefaction wave is not necessarily small. The
stability of the single viscous contact wave is also obtained in
\cite{Qin-Wang} if the contact wave is weak enough. It should be
remarked that the subsonic BL-solution decays exponentially with
respect to $\x=x-\s_-t$, which is good enough to get the desired
estimates. When the boundary value $(\r_-,u_-,\t_-)$ belongs to the
supersonic region, there is no BL-solution. Thus the large-time
behavior of the solution is expected to be same as that of the
Cauchy problem and the stability of the 3-rarefaction waves is also
given in \cite{Huang-Li-Shi}.

In the present paper, we are interested in the stability of wave
patterns to the inflow problem (\ref{(1.1)}), (\ref{(1.3)}) and
$(1.4)_3$ when $(\r_-,u_-,\t_-)$ belongs to the transonic region. In
this case, a new wave structure which contains four waves: the
transonic(or degenerate) BL-solution, 1-rarefaction wave, viscous
2-contact wave and 3-rarefaction wave, occurs. Due to the fact that
the first characteristic speed on the boundary is coincident with
the speed of the moving boundary in the transonic BL-solution case,
the nonlinear waves in the first characteristic field may appear,
which is quite different from the the regime that $(\r_-,u_-,\t_-)$
belongs to the subsonic region in our previous result
\cite{Qin-Wang}, where the waves in the first characteristic field
must be absent. Here we just assume that the 1-rarefaction wave
appear in the first characteristic field. Correspondingly, some new
mathematical difficulties occur due to the degeneracy of the
transonic BL-solution and its interactions with other wave patterns
in the superposition wave. In particular, the transonic boundary
layer solution is attached with 1-rarefaction wave for all time, so
the interaction of these two waves should be carefully treated in
the stability analysis.

Because the system $(\ref{(1.1)})$ we consider  is in one dimension
of the space variable $x$, it is convenient to use the following
Lagrangian coordinate transformation:
$$
(t,x)\Rightarrow \left(t,\int^{(t, x)}_{(0, 0)} \rho(\tau, y)\,
 dy-\rho u(\tau, y)\, d\tau\right).
$$
Thus the system $(\ref{(1.1)})$ can be transformed into the
following moving boundary problem of Navier-Stokes equations in the
Lagrangian coordinates \cite{Matsumura-Nishihara-3}:
\begin{equation}
\begin{cases}
  v_{t}-u _{x} = 0,     \cr
    u_{t} +  p_{x}
    =\m
    \left(\frac{u_{x}}{v}\right) _{x},\qquad \qquad\qquad \qquad\qquad   ~~~ ~~t>0, x > \s_- t, \cr
       \left(\frac{R}{\gamma-1} \theta+\frac{1}{2}u^2\right) _t+    (p u)
       _{x}
    =\k\left( \frac{\theta_{ x}}{v} \right)_{x}+  \m\left(  \frac{ uu_{x}}{v}
    \right)_{x},\cr (v,u, \theta)( 0, x)=(v_0,u_0,\t_0)(x)\rightarrow (v_{+},u_{+},\theta_{+}),~~ ~~ {\rm as} ~~ x\rightarrow
 +\infty, \cr
(v, u,\t)( t, x=\s_- t)=(v_-, u_-, \theta_-),~~u_->0
\end{cases} \label{(1.6)}
\end{equation}
where $\s_-:=-\frac{u_-}{v_-}<0$ is the speed of  the moving
boundary.

In order to fix the moving boundary $x=\s_- t$, we introduce a new
variable $\x=x-\s_- t$. Then we have the half-space problem
\begin{equation}
\begin{cases}
  v_{t}-\s_- v_\x-u _{\x} = 0,     \cr
    u_{t}-\s_- u_\x+  p_\x
    =\m
    \left(\frac{u_{\x}}{v}\right) _{\x}, \qquad\qquad\qquad\qquad\qquad\qquad\qquad\qquad t>0, \x\in\mathbf{R}_+, \cr
       \Big(\frac{R}{\gamma-1} \theta+\frac{1}{2}u^2\Big) _t-\s_- \left(\frac{R}{\gamma-1} \theta+\frac{1}{2}u^2\right) _\x+   (p u) _\x
    =\k\left( \frac{\theta_{  \x}}{v} \right)_\x+  \m\left(  \frac{ uu_\x}{v}
    \right)_\x,\cr(v,u, \theta)(t=0, \x)=(v_{0},u_{0}, \theta_0)(\x)\rightarrow (v_{+},u_{+},\theta_{+}) ~~{\rm as} ~~ \x \rightarrow
 +\infty, \cr (v, u, \theta)(t,\x=0)=(v_-, u_-, \theta_-),~~u_->0.
\end{cases}\label{(1.7)}
\end{equation}
Given the right end state $(v_+,u_+,\t_+)$, we can define the
following wave curves in the phase space $(v,u,\t)$ with $v>0$ and
$\t>0$.

\noindent$\bullet$ Transonic(or degenerate) boundary layer curve:
\begin{equation}
BL(v_+,u_+,\theta_+):=\left\{(v,u,\theta)\bigg|\f{u}{v}=-\s_-=\f{u_-}{v_-},(u,\t)\in
\Sigma(u_+,\t_+)\right\},\label{(1.8)}
\end{equation}
 where $(v_+,u_+,\t_+)\in \G_{trans}^+=\{(u,\t)|u=\sqrt{R\g\t}>0\,\}$
is the transonic region defined in (\ref{(2.3)}) with positive gas
velocity and $\Sigma(u_+,\t_+)$ is the trajectory at the point
$(u_+,\t_+)$ defined in Case II of Lemma 2.1 below.

\noindent$\bullet$ Contact wave curve:
\begin{equation}
 CD(v_+,u_+,\t_+):=\{(v,u,\theta)|u=u_+, p=p_+, v\not\equiv v_+\},\label{(1.9)}
\end{equation}
$\bullet$ $i-$Rarefaction wave curve $(i=1,3)$:
\begin{equation}
 R_i (v_+, u_+, \theta_+):=\left  \{  (v, u, \theta)\bigg{ |}\l_i<\l_{i+} ,~u=u_+-\int^v_{v_+}
 \lambda_i(\eta,
s_+) \,d\eta,~ s(v, \theta)=s_+\right  \} ,\label{(1.10)}
\end{equation}
where $s_+=s(v_+,\t_+)$ and $\l_i=\l_i(v,s)$ is the $i-$th
characteristic speed given in (\ref{(2.1)}).

Our main stability result is, roughly speaking, as follows:

$\bullet$ Assume that $(v_-,u_-,\t_-) \in {\rm
BL\texttt{-}R_1\texttt{-}CD\texttt{-}R_3} (v_+,u_+,\t_+)$, that is,
there exist the unique medium states
$(v_*,u_*,\t_*)\in\Gamma_{trans}^+$, $(v_m,u_m,\t_m)$ and
$(v^*,u^*,\t^*)$, such that $(v_-,u_-,\t_-) \in {\rm
BL}(v_*,u_*,\t_*)$, $(v_*,u_*,\t_*)\in R_1(v_m,u_m,\t_m)$,
$(v_m,u_m,\t_m) \in {\rm CD}(v^*,u^*,\t^*)$ and $(v^*,u^*,\t^*)\in
{\rm R_3} (v_+,u_+,\t_+)$, then the superposition of the four wave
patterns: the transonic (or degenerate) BL-solution, 1-rarefaction
wave, 2-viscous contact wave and 3-rarefaction wave is
time-asymptotically stable provided that the wave strength
$\d=|(v_+-v_-,u_+-u_-,\t_+-\t_-)|$ is suitably small and the
conditions in Theorem 2.1 hold.

This paper is organized as follows. In Section 2, after giving some
preliminaries on boundary layer solution, viscous 2-contact wave,
rarefaction waves and their superposition, we state our main result.
In Section 3, first the wave interaction estimations are shown, then
the desired energy estimates are performed and finally our main
result is proven.

\vskip 2mm \noindent\emph{Notations.} Throughout this paper, several
positive generic constants are denoted by $c,C$ without confusion,
and $C(\cdot)$  stands for some generic constant(s) depending only
on the quantity listed in the parenthesis. For function spaces,
$L^{p}(\mathbf{R}_+), 1\leq p\leq \infty$, denotes the usual
Lebesgue space on $\mathbf{R}_+$. $W^{k,p}(\mathbf{R}_+)$ denotes
the $k^{th}$ order Sobolev space, and if $p=2$, we note
$H^{k}(\mathbf{R}_+):=W^{k,2}(\mathbf{R}_+)$, $\| \cdot \|:=\| \cdot
\|_{L^{2}(\mathbf{R}_+)}$, and $\| \cdot \|_k:=\| \cdot
\|_{H^{k}(\mathbf{R}_+)}$  for simplicity.  The domain
$\mathbf{R}_+$ will be often abbreviated without confusion.

\section{Preliminaries and Main Result}

It is well known that the hyperbolic system (\ref{(1.5)}) has three
characteristic speeds
\begin{eqnarray}
\lambda_1(v,\t)= -\frac{\sqrt{R \gamma\theta}}{v} ,~~~ \lambda_2=0
,~~~ \lambda_3(v,\t)= \frac{\sqrt{R \gamma\theta}}{v} .
\end{eqnarray}
The first and the third characteristic field is genuinely nonlinear,
which may have nonlinear waves, shock wave and rarefaction wave,
while the second  characteristic field is linearly degenerate, where
contact discontinuity may occur.

Let
\begin{eqnarray}
c(v, s)  :=\sqrt{  - v^2 p_v(v, s)  }=
 \sqrt{R \gamma\theta}=:c(v, \theta), \quad   M(v, u, \theta   ) :=\frac{|u
 |}{c(v,\t)
}\label{(2.1)}
\end{eqnarray}
be the sound speed and the Mach number at the state $(v,u,\t)$.
Correspondingly, set
\begin{eqnarray}
c_+ := c(v_+, \theta_+)=\sqrt{ R\gamma\theta_+} , \quad M_+: =
M(v_+, u_+, \theta_+)=\frac{|u_+|}{c_+}
\end{eqnarray}
be the   sound speed and the Mach number at the far field
$\{x=+\infty \}$. We divide the phase space $\{(v,  u, \theta)|\,
v>0, \theta>0\}$ into three parts:
\begin{eqnarray}
\begin{cases}
 ~~\Omega_{sub} :=\left\{(v,  u, \theta)~|~M<1  \,\right\},\cr
 \Gamma_{trans} :=\left\{(v, u, \theta)~| ~M=1 \,\right\},\cr
\Omega_{super} :=\left\{(v, u, \theta)~|~ M>1 \, \right\}.
\end{cases}
\label{(2.3)}
\end{eqnarray}
Call them subsonic, transonic and supersonic region, respectively.
Obviously, if we add the alternative condition $u> 0$ or $u\leq 0$,
then we have six regions $\Omega_{sub}^{\pm}$,
$\Gamma_{trans}^{\pm}$, and $\Omega_{super}^{\pm}$.

\subsection{ Boundary layer solution}

When $(v_-,u_-,\t_-)\in \Omega_{sub}^+\cup \Gamma^+_{trans}$, we
have
\begin{eqnarray}
\l_1(v_-,\t_-)=-\frac{\sqrt{R
\gamma\theta_-}}{v_-}\leq-\frac{u_-}{v_-} = \s_-<0,
\end{eqnarray}
hence a stationary solution $\big(V^b,U^b,\T^b\big)(\x)$ to the
inflow problem (\ref{(1.7)}) is expected
\begin{eqnarray}\label{(2.6)}
\begin{cases}
  -\s_-  V^b_\x-U^b _{\x} = 0,     \cr
     -\s_-   U^b_\x+   P^b_\x
    =\m
    \Big(\frac{U^b_{\x}}{V^b}\Big) _{\x}, \cr
        -\s_-  \left( \f{R}{\g-1}  \T^b+\frac{1}{2}\left(U^{b  }\right)^2\right)  _\x+   \left(P^b U^b\right) _\x
    =\k\Big( \frac{\T^b_{  \x}}{ V^b} \Big)_\x+  \m\Big(  \frac{ U^bU^b_\x}{ V^b}
    \Big)_\x, \cr
 \big(V^b , U^b, \T^b\big)(0)=(v_-, u_-, \theta_-), ~~~\big(V^b, U^b, \T^b\big)
 (+\i)=
 (v_+, u_+, \theta_+),
\end{cases}
\end{eqnarray}
where $ P^{b}:=p\big( V^{b}, \T  ^{b}\big)=\f{R\T^b}{V^b}$. We call
this stationary solution $\big(V^b,U^b,\T^b\big)(\x)$ the boundary
layer solution (simply, BL-solution) to the inflow problem
(\ref{(1.7)}).

 From the
fact that $V^b(\x) >0$ and $ u_->0$,  then
\begin{equation}\label{(2.7)}
u_+ >0 , \qquad  \frac{U^b}{V^b}=
\frac{u_+}{v_+}=\frac{u_-}{v_-}=-\s_-.
\end{equation}
Thus \eqref{(2.6)} is equivalent to \eqref{(2.7)} and the following
ODE system
\begin{eqnarray}
\label{28}\begin{cases}
  \left(U^{b }\right)^\prime= -\frac{\s_-}{\m}V^b\big(U^{b}-u_+\big)
  +\frac{R}{\m}\left(  \T^{b}  - \frac{\t_+}{v_+}V^b\right)   \qquad\quad^\prime=\f{d}{d\x},\cr
\left(\T^{b }\right)^\prime =- \frac{ R\s_-}{\k(\g-1)
}V^b\big(\Theta^b-\t_+\big)+
 \frac{p_+}{\k}V^b\big(U^b-u_+\big)+\frac{\s_-}{2\k}V^b\big(U^b-u_+\big)^2  ,\cr
  \left(U^{b }, \T^{b } \right)(0)=(u_-, \theta_- ), ~~~\left( U^{b }, \T^{b }\right)(+\infty)=  ( u_{+},\theta_{+})
  ,
\end{cases}
\label{(2.8)}
\end{eqnarray}
 where $p_+:=p(v_+,   \t _+) $.

We can compute that the

Now we state the existence results of the BL-solution to
(\ref{(2.8)}) while its proof has been shown in \cite{Qin-Wang}.

\vskip 2mm

\noindent\textbf{Lemma 2.1  (Existence of BL-solution)
\cite{Qin-Wang}}~ \emph{Suppose that $v_\pm>0$, $u_- >0$,
$\theta_\pm>0$ and let  $\d_b :=|( u_+-  u_-,\theta_+ - \theta_-)|$.
If $u_+\leq 0$, then there is no solution to $(\ref{(2.8)})$. If
$u_+> 0$, then there exists a suitably small constant $\delta_0>0$
 such that if $ 0<\delta^b \leq\delta_0$, then the existence and
non-existence of solutions to \eqref{(2.8)} is divided into three
cases according to the location of $(u_{+}, \theta_+)$: }

\emph{Case I : $ (u_{+}, \theta_+)\in \Omega_{sup}^{+}$. Then there
is no solution to (\ref{(2.8)}).}

\emph{Case II : $ (u_{+}, \theta_+)\in \Gamma_{trans}^{+}$. Then
 $(u_{+}, \theta_+)$ is a saddle-knot point to
\eqref{(2.8)}. Precisely,  there exists a unique trajectory $\Sigma$
tangent to the straight line
\begin{eqnarray}
\m u_+(u-u_+)- \k(\g-1)(\t-\theta_+)=0
\end{eqnarray}
at the point $(u_+,\t_+).$ For each $( u_-,\theta_-)\in \Sigma
(u_+,\t_+) $, there exists a unique solution $\big(U^b,
\Theta^b\big) $ satisfying
$$
U^b_\x>0,\qquad \T^b_\x>0,
$$
and
\begin{eqnarray}
\left|\frac{d^n}{d\x^{n}}\big(  U^b    - u_+ , \Theta^b   -
\theta_+\big )\right| = O(1)\frac{\d_b^{n+1}}{(1+\d_b \x)^ {n+1}}
,~~~n=0,1, 2, \dots.\label{(2.12)}
\end{eqnarray}}

\emph{Case III : $ (u_{+}, \theta_+)\in \Omega_{sub}^{+}$.   Then
$(u_{+}, \theta_+)$ is a saddle point to \eqref{(2.8)}. PPrecisely,
there exists a center-stable manifold $\mathcal{M} $ tangent to the
line
$$
(1+a_2c_2u_+)(U^B-u_+) -a_2(\Theta^B-\theta_+)=0
$$
on the opposite directions at the point $(u_{+}, \theta_+)$. Here
$c_2$ is one of the solutions to the equation
$$
 y^2+ \Bigg{(}\frac{M_+^2\gamma-1}{M_+^2R\gamma}-\frac{\m }{\k(\g-1)}\Bigg{)} y-\frac{\m }{M_+^2 R \gamma \k}=0
$$
 and $a_2=-\frac{R}{\m(\l_J^1-\l_J^2) }$ with $\l_J^1>0,~~\l_J^2<0$ are the two eigenvalues of the linearized matrix of ODE \eqref{(2.8)}.
 Only when $(u_-, \theta_-)\in \mathcal{M}(u_+,\t_+)$,
  does there exist   a unique solution $\big(U^b, \Theta^b\big)\subset \mathcal{M}(u_+,\t_+)$   satisfying}
\begin{eqnarray}
\left| \frac{d^n}{d\x^{n}}\big(  U^b    - u_+ , \Theta^b   -
\theta_+ \big)\right|  = O(1) \d_b   e^{-c\x },~~~n=0,1, 2, \dots.
\end{eqnarray}

\subsection{Viscous Contact  Wave} If $(v_-, u_-, \t_-)\in
CD(v_+,u_+,\t_+)$, then the following Riemann problem
\begin{eqnarray}
\begin{cases}
  v_{t}-u _{x} = 0,     \cr
    u_{t} +  p_{x}
    =0,\qquad \qquad\qquad       t>0, x\in\mathbf{R}, \cr
       \left(\frac{R}{\g-1}  \t+\frac{1}{2}u^2\right) _t+    (p u)
       _{x}
    =0,\cr (v,u, \theta)( 0, x)  =
 \begin{cases}
 (v_-, u_-,
\t_-), \quad x<0,\cr(v_{+},u_{+},\theta_{+}),\quad x>0
\end{cases}
\end{cases}
\end{eqnarray}
admits a contact discontinuity solution
$$
(v,u,\t)(t,x)=\left\{
\begin{array}{ll}
(v_-,u_-,\t_-),\qquad & x<0,~ t>0,\\
(v_+,u_+,\t_+),\qquad & x>0,~ t>0.
\end{array}
\right.
$$

From \cite{Huang-Xin-Yang}, the viscous version of the above contact
discontinuity, called viscous contact wave $\big(V^{d}, U^{d},
\T^{d}\big)(t,x)$ can be defined by
\begin{eqnarray}\label{(2.13)}
\begin{cases}
   V^{d}(t, x)
    =  \frac{R\T^{\rm sim}(t, x)}{p_+ }  , \cr
    U^{d}(t, x)
    = u_++    \frac{(\g-1)\k\T^{\rm sim}_x(t, x)}{\g\T^{\rm sim}(t, x)
    }, \cr \T^{d}(t, x)
    = \T^{\rm sim}\Big(\frac{x}{\sqrt{1+t}}\Big)+R \Big(\mu-\frac{(\g-1)\k}{R  \g}\Big)\T^{\rm sim}_t
\end{cases}
\end{eqnarray}
where  $\T^{\rm sim}\left(\frac{x}{\sqrt{1+t}}\right)$ is the unique
self-similar solution to the following nonlinear diffusion equation
\begin{eqnarray}
\begin{cases}
 \T_t
    = \frac{(\g-1)\k p_+}{   R^2 \g}\left(\frac{\T_x}{\T}\right)_x, \cr
    \T(t,\pm\infty)=\t_{\pm}.
\end{cases}
\end{eqnarray}
Note that $\x=x-\s_- t$, we have the following Lemma:

\noindent\textbf{Lemma 2.2. }\cite{Huang-Xin-Yang}  \emph{The
viscous contact wave  $\big(V^{d}, U^{d},
\T^{d}\big)(t,x),~(x=\x+\s_-t)$ defined in \eqref{(2.13)} satisfies
\newcounter{cd}
\begin{list}
{\upshape  \roman{cd}) \;} {\setlength{\parsep}{\parskip}
 \setlength{\itemsep}{0ex plus0.1ex}
 \setlength{\rightmargin}{1em}
 \setlength{\leftmargin}{3em}
 \setlength{\labelwidth}{4em}
 \setlength{\labelsep}{0.1em}
 \usecounter{cd}\setcounter{cd}{0}
}
\item  $\partial_\x^n\big(\T^{d}-\t_\pm \big)
   =O(1)\d_d(1+t)^{-\frac{n}{2}}\exp \left( -\frac{C_d (\x+\s_-  t)^2}{1+t} \right), \quad n=0,
   1,2,\cdots $;
\item $U^{d}_\x(t, \x)
  =O(1)\d_d(1+t)^{-1}\exp \left( -\frac{C_d (\x+\s_-  t)^2}{1+t}
   \right);$
 \item  $ \big(V^{d}, U^{d}, \T^{d}\big)(t,\x= 0)-(v_-, u_-, \t_-)=O(1) \d_d e^{-ct} $.
\end{list}
where $\d_d=|\t_{+}-\t_{-}|$ is the amplitude of the viscous contact
wave and $C_d,c>0 $ are constants. }

Then
 the viscous contact wave  $\big(V^{d}, U^{d}, \T^{d}\big)$ defined in \eqref{(2.13)}
 satisfies the system
\begin{eqnarray}
\begin{cases}
  V^{d}_t-\s_-   V^{d}_\x-U^{d} _{\x} = 0,     \cr
    U^{d}_t-\s_-   U^{d}_\x+ P^{d}_\x
    =\m\Big( \frac{U^{d}_{  \x}}{ V^{d}} \Big)_\x,\qquad \qquad   ~ ~~~ t>0,  \x\in\mathbf{R}_+, \cr
        \frac{R}{\g-1}  \big( \T^{d}_t-\s_-  \T^{d}  _\x\big)+    P^{d} U^{d}  _\x
    =\k\bigg( \frac{\T^{d}_{  \x}}{ V^{d}} \bigg)_\x+\m\frac{(U^{d}_{  \x})^2}{ V^{d}}
    +H^d
\end{cases}
\end{eqnarray}
where $ P^{d}:=p\big( V^{d}, \T^{d}\big)$ and
$$
H^d=O(1)\d_d(1+t)^{-2}\exp \left( -\frac{C_d(\x+\s_-
t)^2}{1+t} \right)
$$
due to Lemma 2.2.

\subsection{ Rarefaction waves}

It is well known  that if $(v_-, u_-, \t_-)\in R_i (v_+, u_+,
\theta_+),~(i=1,3)$, then there exist a $i-$rarefaction wave
$(v^{r_i}, u^{r_i}, \t^{r_i})(x/t)$ which is the global  weak
solution to the following Riemann problem
\begin{eqnarray}
\begin{cases}v  _{t}- u  _{x} = 0,     \cr
    u_{t}+   p _x
    =0, \quad\quad\quad\quad\quad\quad  \,t>0,  x\in\mathbf{R}, \cr
     \left( \frac{R}{\g-1}  \t+\frac{1}{2}u^2\right) _t+ (   pu )_x
    =0,  \cr
 (v,u, \theta)( 0, x)=
 \begin{cases}
 (v_-, u_-,
\t_-), \quad x<0,\cr(v_{+},u_{+},\theta_{+}),\quad x>0 .
\end{cases}
\end{cases}
\end{eqnarray}
Consider the following Burgers equation
\begin{eqnarray}
\begin{cases}
 w_{t}+ww_{x}=0,  \quad \,t>0,  x\in\mathbf{R}, \cr
 w_0( x):= w( 0,x)=
\begin{cases}
w_-, ~\quad \quad\quad\quad\quad\quad \quad\quad   x<0 ,\cr\di
w_-+C_q (w_+-w_-) \int^{ x }_0y^qe^{-y}\,dy,~ x\geq 0 .
\end{cases}
\end{cases}
\end{eqnarray}
Here $ q\geq 14$ is a constant to be determined, and $C_q$ is a
constant such that $\di C_q\int^{+\infty}_0y^qe^{-y} dy=1 $. If
$w_-<w_+,$ then the solution to the above Burgers equation can be
expressed by
\begin{eqnarray}
w(t,x)=w_0(x_0(t, x)),\quad\quad x=x_0(t, x)+w_0(x_0(t, x))t.
\end{eqnarray}
Moreover, we have

$\bullet$ $w(t,x)=w_-$, if $x\leq w_-t$.

$\bullet$ For any positive constant $\s_0>0$ and for $x\geq0$
\begin{eqnarray}\label{(2.19)}
|w(t, x)-w_+|&=&|w_0(x_0(t,x))-w_+|\cr\cr &= &C_q (w_+-w_-) \int_{
x_0(t,x)}^{+\infty} y^qe^{-y}\,dy\cr\cr &= &C_q (w_+-w_-)
\int_{x-w_0( x_0(t,x))t}^{+\infty} y^qe^{-y}\,dy \cr\cr&\leq &C_q
(w_+-w_-) \int_{x-w_+t}^{+\infty} y^qe^{-y}\,dy\cr\cr &\leq  &C_q
(w_+-w_-)e ^{ -\s_0  t },
  \qquad \text{if}~~x\geq (2\s_0+w_+) t.
\end{eqnarray}
Note that the estimation in \eqref{(2.19)} play an important role in
the wave interaction estimates, which is motivated by
\cite{Liu-Matsumura-Nishihara} and \cite{Matsumura-Nishihara-1} .

Now  the $i-$rarefaction wave $(V^{r_i}, U^{r_i}, \T^{r_i})
(t,x)~(i=1,3)$ to the inflow problem \eqref{(1.7)} can be defined by
\begin{eqnarray}\label{(2.20)}
\begin{cases}
 \l_i(V^{r_i},\T^{r_i})  (t,x)  =w(1+t,x+\s_-) ,\cr
 s(V^{r_i}, \T^{r_i})(t,x)=s_+=s(v_+,\theta_+),  \cr \di
U^{r_i}(t,x)=u_+-\int^{V^{r_i}(t,x)}_{ v_+} \l_i(\eta,s_+) d\eta.
\end{cases}
\end{eqnarray}
Then the $i-$rarefaction wave $(V^{r_i}, U^{r_i}, \T^{r_i})(t,x), ~
(i=1, 3)$ defined in \eqref{(2.20)} satisfies the system
\begin{eqnarray}\label{(2.21)}
\begin{cases}
  V^{r_i}_t-\s_-  V^{r_i}_\x-U^{r_i} _{\x} = 0,     \cr
    U^{r_i}_t- \s_-  U^{r_i}_\x+   P^{r_i}_\x
    =0,\cr
        \left[ \frac{R}{\g-1}   \T^{r_i}+\frac{1}{2}(U^{{r_i}  })^2\right]  _t
        -\s_- \left[\frac{R}{\g-1}   \T^{r_i}+\frac{1}{2}(U^{{r_i}  } )^2\right]  _\x+   (P^{r_i} U^{r_i}) _\x
    =0, \cr
 (V^{r_i} , U^{r_i}, \T^{r_i})(t, \x=0)=(v_-, u_-, \theta_-),\cr (V^{r_i} , U^{r_i},
 \T^{r_i}) (t, \x)\rightarrow (v_{+},u_{+},\theta_{+}) ~~\text{ as} ~~ \x \rightarrow
 +\infty
\end{cases}
\end{eqnarray}
where $ P^{r_i}:=p( V^{r_i}, \T  ^{r_i})$.

 \vskip 2mm
\noindent\textbf{Lemma 2.3} \emph{ $i-$rarefaction wave $(V^{r_i},
U^{r_i}, \T^{r_i})(t,\x), ~ (i=1, 3)$ defined in \eqref{(2.20)}
satisfies}
\newcounter{w}
\begin{list}
{\upshape  \roman{w}) \;} {\setlength{\parsep}{\parskip}
 \setlength{\itemsep}{0ex plus0.1ex}
 \setlength{\rightmargin}{1em}
 \setlength{\leftmargin}{3em}
 \setlength{\labelwidth}{4em}
 \setlength{\labelsep}{0.1em}
 \usecounter{w}\setcounter{w}{0}
}
\item \emph{  $U^{r_i}_{\x}(t,\x)>0,~~(|V^{r_i}_\x|, |\T_\x^{r_i}|)\leq C U^{r_i}_{\x}$};
\item  \emph{For any $p$} ($1\leq p\leq \infty$), \emph{there exists a constant
$C_{pq}$ such that
\begin{eqnarray}\nonumber
&& \| (V^{r_i}_\x, U^{r_i}_\x, \T^{r_i}_\x)(t)\|_{L^p}\leq C_{p}
\min\big\{\d_{r_i}, \d_{r_i}^{1/p}(1+t)^{-1+1/p}\big\} , \cr&&\|
(V^{r_i}_{\x\x}, U^{r_i}_{\x\x}, \T^{r_i}_{\x\x})(t)\|_{L^p}\leq
C_{p} \min\big\{\d_{r_i},\d_{r_i}^{1/p+1/q} (1+t)^{-1+1 /q } \big\};
 \end{eqnarray}
  \item  For $\forall\,\s_0>0$, if\, $   \x\geq \left[-\s_- +\l_1(v_+,\t_+) +2\s_0\right](1+t)
 $, then
   $\Big|\partial_\x^n\big\{(V^{r_1}, U^{r_1}, \T^{r_1})(t,\x) -(v_+, u_+, \t_+)\big\}\Big|
\leq C \d_{r_1} e^{-\s_0 t},~n=0,1,2,\cdots;$
\item  For $  \x\leq \left[-\s_-+\l_3(v_-,\t_-)\right](1+t) $,
   $ (V^{r_3}, U^{r_3}, \T^{r_3}) -(v_-, u_-, \t_-) \equiv 0;$
\item   $ \lim\limits_{t\rightarrow
\infty}\sup\limits_{\x\in\mathbf{R}_+}\big| (V^{r_i}, U^{r_i},
\T^{r_i})( t, \x)-(v^{r_i},u^{r_i}, \theta^{r_i})
\big(\frac{\x}{1+t}\big)\big| =0$.}
\end{list}

\noindent{\bf Remark:} The statement ${\rm iii)}$ is a direct
consequence of the \eqref{(2.19)}.

\subsection{Superposition of transonic BL-solution, 1-rarefaction
wave, 2-viscous contact  wave and 3-rarefaction wave} In this
subsection, we consider the case that $ (v_-,u_-,\t_-)  \in
BL$-$R_1$-$CD$-$R_3(v_+,u_+,\t_+)$, that is, there exist uniquely
three medium states $(v_*, u_*,\t_*)\in \G_{trans}^+$, $(v_m,
u_m,\t_m)$ and $(v^*, u^*,\t^*)$ such that $(v_*, u_*,\t_*)\in
BL(v_-, u_-,\t_-)$, $(v_*, u_*,\t_*)\in R_1(v_m,u_m,\t_m)$, $(v_m,
u_m,\t_m)\in CD(v^*, u^*,\t^*)$  and $(v^*, u^*,\t^*)\in
R_3(v_+,u_+,\t_+)$. In fact, three medium states $(v_*, u_*,\t_*)\in
\G_{trans}^+$, $(v_m, u_m,\t_m)$ and $(v^*, u^*,\t^*)$ can be
expressed explicitly and uniquely by the following nine equations
\begin{eqnarray}\label{(2.22)}
\begin{cases}
 \di \frac{u_-}{v_-}= \frac{u_*}{v_*}, \quad  u_*=\sqrt{R\g \t_*}, \quad
 (u_-,\t_-)\in \Sigma (u_*, \t_*),  \cr
   \di u_*=u_m-\int_{v_*}^{v_m} \sqrt{R\gamma v_+^{ \g-1} \t_+}~\eta^{
-\frac{ \gamma+1}{2}} \,d\eta, \quad ~ v_* ^{  \g-1}  \t_* = v_m ^{
\g-1} \t_m, \cr
       \di u_m=u^*, \quad  \frac{\t_m}{v_m}= \frac{\t^*}{v^*},  \cr \di u^*=u_++\int_{v^*}^{v_+} \sqrt{R\gamma v_+^{ \g-1} \t_+}~\eta^{
-\frac{ \gamma+1}{2}} \,d\eta, \quad ~ v^{*   \g-1}  \t^* =v_+^{
\g-1} \theta_+.
\end{cases}
\end{eqnarray}
Define the superposition wave $(V,U,\T)(t,\x)$ by
\begin{eqnarray}
 \left(\begin{array}{cc} V\\ U \\
  \T
\end{array}
\right)(t, \x)= \left(\begin{array}{cc} V^b+ V^{r_1}+ V^{d}+ V^{r_3}\\ U^b + U^{r_1}+ U^{d}+ U^{r_3} \\
\T^b+ \T^{r_1}+ \T^{d}+ \T^{r_3}
\end{array}
\right)(t, \x) -\left(\begin{array}{cc} v_*+v_m+v^*\\ u_*+u_m+u^*\\
\t_*+\t_m+\t^*
\end{array}
\right) \label{(2.28)}
\end{eqnarray}
where $(V^b, U^b, \Theta^b )(\x)$ is the transonic BL-solution
defined in Case II of Lemma 2.1  with the right state $(v_+, u_+,
\t_+)$ replaced by $(v_*, u_*, \theta_* )$, $(V^{r_1}, U^{r_1},
\Theta^{r_1} )(t,\x)$ is the 1-rarefaction wave defined in
(\ref{(2.20)}) with the states $(v_-, u_-, \t_-)$ and $(v_+, u_+,
\t_+)$  replaced by $(v_*, u_*, \theta_* )$ and $(v_m, u_m, \t_m)$
respectively, $(V^{d}, U^{d}, \Theta^{d} )(t,\x)$ is the viscous
contact wave defined in (\ref{(2.13)}) with the states $(v_-, u_-,
\t_-)$ and $(v_+, u_+, \t_+)$ replaced by $(v_m, u_m, \theta_m )$
and $(v^*, u^*, \theta^* )$, respectively, and $(V^{r_3}, U^{r_3},
\Theta^{r_3} )(t,\x)$ is the 3-rarefaction wave defined in
(\ref{(2.20)}) with the left state $(v_-, u_-, \t_-)$ replaced by
$(v^*, u^*, \theta^* )$.

Now we state the main result of the paper as follows.

\

\noindent \textbf{Theorem 2.1 (Stability of superposition of four
waves)}~\emph{Assume that $(v_-,u_-, \theta_-)\in BL
$-$R_1$-$CD$-$R_3(v_+,u_+,\t_+)$.  Let $(V,U,\T)(t,\x)$ be the
superposition of the transonic BL-solution, 1-rarefaction wave,
viscous 2-contact wave and 3-rarefaction wave defined in
(\ref{(2.28)}). Then there exists a small positive constant $\d_0$
such that if the   initial values and the wave strength
$\d=|(v_+-v_-, u_+-u_-, \t_+-\t_-)|$ satisfy
\begin{eqnarray}
 \d +\|(v_{0}-V_0,   u_{0}- U_0, \theta_0-\T_0)\|_1  \leq \d_0.
\end{eqnarray}
the inflow problem $(\ref{(1.7)})$ has a unique global-in-time
solution $(v, u, \theta)(t,\x)$ satisfying}
\begin{eqnarray}\begin{cases}
 (v-V , u-U  , \t - \T )  ( t,\x) \in C\big([0,\infty);H^{1}(\mathbf{R}^+)\big),\cr
 (v-V ) _{\x} ( t, \x )\in
 L^{2}\big(0,\infty;L^{2}(\mathbf{R}^+)\big),\cr
  ( u-U  , \t - \T)  _\x( t, \x )\in
  L^{2}\big(0,\infty;H^{1}(\mathbf{R}^+)\big).
\end{cases}\end{eqnarray}
\emph{ Furthermore,}
\begin{eqnarray}
\lim_{t\rightarrow \infty}\sup_{\x\in
\mathbf{R}_+}|(v-V,u-U,\t-\T)(t, \x)| =0.
\end{eqnarray}

\

\noindent \textbf{Remark.}  In Theorem 2.1, we assume that
$\d=|(v_+-v_-, u_+-u_-, \t_+-\t_-)|$ is suitably small. This
assumption is equivalent to the one that the amplitudes of the four
waves are all suitably small. In fact, from the relations in
\eqref{(2.22)} and the facts $U^{b}_\x>0$, $U^{r_1}_\x>0$,
$U^{r_3}_\x>0$, we have
\begin{eqnarray}
\begin{cases}
  | v_*-v_-|+| \t_*-\t_-|  =O(1)( u_*-u_-) ,
  \cr    | v_m-v_*|+| \t_m-\t_*|  =O(1)(u_m-u_*),
 \cr  |v_+- v^* |+|\t_+- \t^* |  =O(1)(u_+- u^* ).
\end{cases}
\end{eqnarray}
Thus $\d_b=O(1)( u_*-u_-), \d_{r_1}=O(1)(u_m-u_*)$,
$\d_{r_3}=O(1)(u_+-u^*)$. Due to $u_m=u^*$ by the contact
discontinuity curve, we have if $\d$ is small, then $\d_b, \d_{r_1}$
and $\d_{r_3}$ are all small. Furthermore, we have
$\d_{d}=|\t^*-\t_m|\leq \d_b+\d_{r_1}+\d_{r_3} +\d$ is small.

\section{Stability Analysis}

\subsection{Wave interaction estimates}

Recalling the definition of the superposition wave $(V,U,\T)(t,\x)$
defined in (\ref{(2.28)}),   we have
\begin{eqnarray}\label{(3.1)}
\begin{cases}
  V_t-\s_-  V_\x-U _{\x} = 0,     \cr
    U_t-\s_-   U_\x+   P_\x
    = \m \Big(\frac{U_{\x}}{V}\Big) _{\x}+G,\qquad \qquad\qquad   ~~~~ ~    t>0,  \x\in\mathbf{R}_+, \cr
        \f{R}{\g-1} \left( \T _t-\s_-    \T  _\x\right)+   P U _\x
    =\k\Big( \frac{\T_{  \x}}{ V} \Big)_\x+  \m    \frac{  (U_\x)^2}{ V}
     +H,  \cr (V, U, \T)  (t, \x=0)=(v_-, u_-,
     \t_-)+\left(V^{d},U^{d},\T^{d}\right)(t,\x=
    0)-(v_m, u_m, \t_m).
\end{cases}
\end{eqnarray}
 where $ P :=p( V , \T  )$ and
\begin{eqnarray}
\begin{cases}
  G= \big( P-P^b-P^{r_1}-P^{d}-P^{r_3}\big )_\x- \m \bigg(\frac{U_{\x}}{V}-\frac{U^b_{\x}}{V^b}-\frac{U^d_{\x}}{V^d}\bigg) _{\x}
   =:G_1+ G_2, \cr
       H =( PU_\x-P^bU^b_\x -P^{r_1} U^{r_1}_\x-P^{d}U^{d}_\x
       -P^{r_3} U^{r_3}_\x)\cr~~~~~~
       -\bigg[\k\bigg(\frac{\T_{  \x}}{ V} - \frac{\T^b_{  \x}}{ V^b}
       -\frac{\T^{d}_{  \x}}{ V^{d}}\bigg)_\x+ \m  \Bigg( \frac{  (U_\x)^2}{ V}- \frac{  \big(U^{b }_\x\big)^2}{ V^b}
    -   \frac{  \big(U^{d }_\x\big)^2}{ V^d}\Bigg)-H^d\bigg]=:H_1+ H_2.
\end{cases}
\end{eqnarray}

To control the interaction terms coming from different wave
patterns, we give the following lemma which  will be critical in the
energy estimate in Subsection 3.3.

\vskip 2mm
 \noindent\textbf{Lemma 3.1 (Wave interaction estimates)} \emph{ }
\begin{eqnarray}\label{3.3}
\begin{cases}
  \di \int_{\mathbf{R}_+}
  \big|V^b_\x  \big(  V^{r_1}-v_*\big)\big|+\big|V^{r_1}_\x \big(V^b-v_*\big)\big|
 \,d\x =O(1) \d^{1/8} ( 1+t)^{ -13/16}   ,\cr
\di  \int_{\mathbf{R}_+} \big|V^b_\x  \big(V^{d}-v_m\big)
\big|+\big|V^{d}_\x \big (V^b-v_*\big) \big|\, d\x  =O(1) \d ( 1+
t)^{-1},
 \cr \di  \int_{\mathbf{R}_+}
  \big|V^b_\x  \big(  V^{r_3}-v^*\big)\big|+\big|V^{r_3}_\x \big(V^b-v_*\big)\big|
 \,d\x =O(1)  \d^{1/8} ( 1+t)^{ -7/8} ,
 \cr\di \int_{\mathbf{R}_+} \big|V^{d}_\x \big (V^{r_1}-v_m\big)
\big|+\big|V^{r_1}_\x  \big(V^{d}-v_m\big) \big|\, d\x  =O(1) \d
e^{- c t },
 \cr
\di \int_{\mathbf{R}_+} \big|V^{d}_\x \big(V^{r_3}-v^*\big)\big
|+\big|V^{r_3}_\x \big(V^{d}-v^*\big) \big|\, d\x=O(1)\d e^{- c t }
, \cr\di  \int_{\mathbf{R}_+} \big|V^{r_1}_\x
\big(V^{{r_3}}-v^*\big)\big|+ \big|V^{r_3}_\x \big(  V^{{r_1}}-v_m
\big)\big|\, d\x =O(1)\d e^{-c  t} ,
\end{cases}
\end{eqnarray}
\begin{eqnarray}\label{(3.4)}
\begin{cases}
 \di \int_{\mathbf{R}_+} \big|V^b_\x   V^d_\x \big|\,
d\x =O(1)\d  ( 1+
 t)^{-2}  , \qquad \int_{\mathbf{R}_+}
  \big|V^b_\x   V^{r_1}_\x\big|
 \,d\x =O(1) \d( 1+t)^{- 1} , \cr \di  \int_{\mathbf{R}_+}
  \big|V^b_\x     V^{r_3}_\x\big|
 \,d\x=O(1) \d ( 1+  t)^{-1}, \,~\quad \int_{\mathbf{R}_+}\big |V^d_\x   V^{r_1}_\x\big |\, d\x
 =O(1)\d
  e^{-c  t},
 \cr\di  \int_{\mathbf{R}_+} \big|V^d_\x
 V^{r_3}_\x \big|\, d\x=O(1)\d e^{-ct} , \quad~~\qquad
\int_{\mathbf{R}_+} \big|V^{r_1}_\x  V^{r_3}_\x
  \big|\, d\x =O(1) \d e^{-c t},
\end{cases}
\end{eqnarray}

\noindent \textbf{Proof.}  First we prove $(\ref{3.3})_1$, that is

$\bullet$ Interaction of transonic boundary layer solution and
1-rarefaction wave:

Since $V^{r_1}_\x\geq 0$ and $V^b_\x \geq 0$, we have  $
V^{r_1}-v_*\geq0$ and $ v_* -V^{b}\geq0$. Thus we have
\begin{eqnarray}\label{(3.6)}
&& \int_{\mathbf{R}_+}
  \big|V^b_\x  \big(  V^{r_1}-v_*\big)\big|+\big|V^{r_1}_\x \big(V^b-v_*\big)\big|
 \,d\x\cr &= & 2\left \{ \int_{0}^{ \left[-\s_-+\f{\l_1(v_m,\t_m)}{2} \right]( 1+t)   }+
 \int _{ \left[ -\s_-+\f{\l_1(v_m,\t_m)}{2} \right](1+t) }^{+\infty} \right\}
  V^{r_1}_\x \big( v_*-V^b\big) \,d\x\cr
  &:=&J_1+J_2.
 \end{eqnarray}
Note that
$$
\begin{array}{ll}
\di
-\s_-+\f{\l_1(v_m,\t_m)}{2}\\
\di=\f{u_-}{v_-}+\f{\l_1(v_m,\t_m)}{2} =\f{u_*}{v_*}+\f{\l_1(v_m,\t_m)}{2}\\
=\f{\sqrt{R\g\t_*}}{v_*}+\f{\l_1(v_m,\t_m)}{2}=-\l_1(v_*,\t_*)+\f{\l_1(v_m,\t_m)}{2}\\
=\left[\l_1(v_m,\t_m)-\l_1(v_*,\t_*)\right]-\f{\l_1(v_m,\t_m)}{2}\\
\geq-\f{\l_1(v_m,\t_m)}{2}>0.
\end{array}
$$
Now we can compute that
\begin{eqnarray}\label{(3.7)}
J_1&=&  \int_{0}^{ \left[ -\s_-+\f{\l_1(v_m,\t_m)}{2} \right]( 1+t)
} V^{r_1}_\x \big( v_*-V^b\big) \,d\x
  \cr &= &   O(1) \|V_\x^{r_1}(t)\|_{L^\infty}
   \int_{0}^{ \left[ -\s_-+\f{\l_1(v_m,\t_m)}{2} \right]( 1+t)   }\frac{\d_b  }{ 1+\d_b\x
}\, d\x \cr &=& O(1)\d_{r_1}^{\frac{1}{8}}  ( 1+t)^{
-\frac{7}{8}}\ln ( 1+\d_b t) \cr&=& O(1)\d_{r_1}^{\frac{1}{8}} (
1+t)^{ -\frac{13}{16}} ,
 \end{eqnarray}
 and
 \begin{eqnarray}\label{(3.8)}
J_2&=&   \int _{ \left[ -\s_-+\f{\l_1(v_m,\t_m)}{2}\right] ( 1+t)
}^\infty
  V^{r_1}_\x \big( v_*-V^b\big) \,d\x
  \cr &= &   O(1)\d_b(v_m-V^{r_1} (t, \x))|_{\x=\left[-\s_-+\f{\l_1(v_m,\t_m)}{2}\right] (
  1+t)}\cr
&=&O(1)\d_be^{-\s_0 t}.
 \end{eqnarray}
due to the statement ${\rm iii)}$ in Lemma 2.3 by taking
$\s_0=-\f{\l_1(v_m,\t_m)}{2}>0$. So the combination of \eqref{(3.7)}
and \eqref{(3.8)} gives $\eqref{3.3}_1$.

Then we prove  $(\ref{3.3})_2$:

$\bullet$ Interaction of transonic boundary layer solution and
viscous 2-contact wave:
$$
\begin{array}{ll}
\di\int_{\mathbf{R}_+} \big|V^b_\x  \big(V^{d}-v_m\big)
\big|+\big|V^{d}_\x \big (V^b-v_*\big) \big|\,
d\x\\
\di=\left\{\int_{0}^{-\frac{\s_- t }{2} }+\int _{-\frac{\s_- t }{2}
}^{+\infty} \right\} \big|V^b_\x \big(V^{d}-v_m\big)
\big|+\big|V^{d}_\x \big (V^b-v_*\big) \big|\, d\x\\
\di :=J_3+J_4.
\end{array}
$$
We calculate
\begin{eqnarray}\label{(3.9)}
J_3&=& \int_{0}^{-\frac{\s_- t }{2} }  \big |V^b_\x  \big
(V^d-v_m\big) \big|+\big|V^d_\x \big(V^b-v_*\big)\big |\, d\x
    \cr &= &O(1) \d_d \int_{0}^{-\frac{\s_- t }{2} }  \exp \left( -\frac{C_d (\x+\s_-  t)^2}{1+t} \right) d\x\cr
    &=& O(1) \d_d     e^{-ct}.
\end{eqnarray}
Also, we have
$$
\begin{array}{l}
\di J_4=\int _{-\frac{\s_- t }{2} }^{+\infty} \big|V^b_\x
\big(V^{d}-v_m\big) \big|+\big|V^{d}_\x \big (V^b-v_*\big) \big|\,
d\x\\\di\quad :=J_4^1+J_4^2.
\end{array}
$$
We can estimate
\begin{eqnarray}\label{(3.10)}
J_4^1&=& \int _{-\frac{\s_- t }{2} }^\infty    \big|V^b_\x
\big(V^d-v_m\big) \big|\, d\x\cr
    &= & O(1) \d_d \d_b^2 ( 1+ \d_b t)^{-2} \int _{-\frac{\s_- t }{2} }^\infty
    \exp \left( -\frac{C_d (\x+\s_-  t)^2}{1+t} \right)  d\x
    \cr &= &   O(1) \d_d   ( 1+  t)^{-3/2} \int _{-\infty }^\infty
    \exp \left( -C_d  \eta^2 \right)  d\eta
   \cr &= &   O(1) \d_d   ( 1+  t)^{-3/2} ,
\end{eqnarray}
and
\begin{eqnarray}\label{(3.11)}
J_4^2&=& \int _{-\frac{\s_- t }{2} }^\infty    \big|V^d_\x
\big(V^b-v_*\big) \big|\, d\x\cr
 &= & O(1) \d_d \d_b ( 1+ \d_b t)^{-1} ( 1+  t)^{-1/2}\int _{-\frac{\s_- t }{2} }^\infty
    \exp \left( -\frac{C_d (\x+\s_-  t)^2}{1+t} \right)  d\x
    \cr &= &   O(1) \d_d   ( 1+  t)^{-1}.
\end{eqnarray}
Thus we proved $\eqref{3.3}_2.$

Now we compute $(\ref{3.3})_3$:

$\bullet$ Interaction of transonic boundary layer solution and
3-rarefaction wave:
\begin{eqnarray}
    &&\int_{\mathbf{R}_+}
  \big|V^b_\x  \big(  V^{r_3}-v^*\big)\big|+\big|V^{r_3}_\x \big(V^b-v_*\big)\big|
 \,d\x  \cr &= & \int _{ \left[
-\s_- + \l_3(v^*,\t^*) \right] ( 1+t) }^\infty
  V^b_\x \big (  v^*- V^{r_3}\big) + V^{r_3}_\x \big(V^b-v_*\big)
 \,d\x
 \cr &= & O(1) \d_b (1+\d_b t )^{-1}  \cr&=& O(1) \min\big\{\d ,(1+  t )^{-1}\big\}\cr
 &=&O(1)\d^{\f18}(1+t)^{-\f78}.
\end{eqnarray}
where in the first equality we have used the fact ${\rm iv)}$ in
Lemma 2.3.

Then we verify $(\ref{3.3})_4$:

$\bullet$ Interaction of 1-rarefaction wave and viscous 2-contact
wave:

First we have
$$\begin{array}{ll}
\di \int_{\mathbf{R}_+} \big|V^d_\x ( v_m-V^{r_1}) \big |\, d\x\\
\di =\left\{\int_{0}^{\left[ -\s_-+\frac{\l_1(v_m,\t_m)}{2} \right]
( 1+t)} +\int _{\left[ -\s_-+\frac{\l_1(v_m,\t_m)}{2} \right](1+t)
}^{+\infty}\right\}
\big|V^d_\x \big(v_m-V^{r_1} \big)    \big|\, d\x\\
\di:=J_5+J_6.
\end{array}
$$
Then we can compute
\begin{eqnarray}
 J_5&=& \int_{0}^{\left[ -\s_-+\frac{\l_1(v_m,\t_m)}{2} \right] (1+t)}  \big|V^d_\x \big(v_m-V^{r_1} \big)    \big|\, d\x
\cr&=& O(1) \d_d  (1+t)^{- \frac{1}{2}} \int_{0}^{\left[
-\s_-+\frac{\l_1(v_m,\t_m)}{2} \right](1+t)}  \exp \left( -
\frac{C_d(\x+\s_- t)^2}{1+t}\right) d\x\cr &=&O(1) \d_d e^{-ct} ,
\end{eqnarray}
and
\begin{eqnarray}
J_6&=& \int _{\left[ -\s_-+\frac{\l_1(v_m,\t_m)}{2} \right] (1+t)
}^\infty   \big|V^d_\x \big(v_m-V^{r_1} \big)    \big|\, d\x \cr&=&
O(1)
  \d_d  \sup _{\x\geq
\left[ -\s_-+\frac{\l_1(v_m,\t_m)}{2} \right] ( 1+t) }\big(
v_m-V^{r_1}(t,\x)\big)=
   O(1)   \d_d \, e^{- c  t } .
\end{eqnarray}
Similarly, we can estimate the interaction term
 \begin{eqnarray}
  \int_{\mathbf{R}_+} \big|V^{r_1}_\x  \big(V^d-v_m\big)
\big|\, d\x  =   O(1) \d_d \, e^{- c t } .
\end{eqnarray}
So $(\ref{3.3})_4$ is verified.

For $(\ref{3.3})_5$, that is

$\bullet$ Interaction of 3-rarefaction wave and viscous 2-contact
wave,  which can be done similarly as $(\ref{3.3})_4$, we omit the
details for simplicity.

Finally, we prove $(\ref{3.3})_6$:

$\bullet$ Interaction of 1-rarefaction wave and  3-rarefaction wave:

Since $V^{r_1}_\x\geq 0$, $V^{r_3}_\x\leq 0$ and the facts ${\rm
iii)}$ and ${\rm iv)}$ in Lemma 2.3, one has
 \begin{eqnarray}
 &&\int_{\mathbf{R}_+}\big |V^{r_1}_\x \big(V^{{r_3}}-v^*\big)\big|+\big |V^{r_3}_\x
 \big(V^{{r_1}}-v_m
\big)\big|\, d\x \cr\cr&=& 2 \int _{\left[-\s_-
+\l_3(v^*,\t^*)\right](1+t)}^{+\infty} V^{r_1}_\x
\big(v^*-V^{{r_3}}\big) \, d\x\cr &=& O(1) \d_{r_1} e^{-c t} =
O(1) \d e^{-c t}.
\end{eqnarray}
Thus we justified \eqref{3.3}. The proof of \eqref{(3.4)} can be
done similarly, but the decay rates with respect to the time $t$ may
be higher. Therefore, we complete the proof of the wave interaction
estimates in Lemma 3.1. $\blacksquare$

With the wave interaction estimation Lemma 3.1 in hand, we have the
following Lemma:

\noindent{\bf Lemma 3.2.} \begin{eqnarray} &&\di  \| G(t) \|_{L^1} +
\| H(t)
 \|_{L^1} =O(1)  \d^{ \frac{1}{8}} ( 1+t)^{ -\frac{13}{16}},\cr&&\di  \| G (t)\| + \| H(t)
 \| =O(1)  \d ( 1+t)^{ -1}   .
\end{eqnarray}

\noindent \textbf{Proof.} ~ We can compute
 \begin{eqnarray}
G_1&=&\big|\big( P-P^b-P^{r_1}-P^{d}-P^{r_3}\big )_\x\big|  \cr
&=&O(1) \big|V^{b} _\x\big| \big(|V^{r_1}-v_* |+\big|V^d-v_m\big|+
      |V^{r_3}-v^*  |\big) \cr
&&+O(1)\big|V^{d} _\x\big| \big(\big|V^b-v_*\big|+ |V^{r_1}-v_m |+
|V^{r_3}-v^*|\big)\cr
&&+O(1)\big|V^{r_1}_\x\big|\big(\big|V^b-v_*\big|+\big|V^{d}-v_m\big|+|V^{r_3}-v^*|\big)\cr
&&+O(1)\big|V^{r_3}_\x\big|\big(\big|V^b-v_*\big|+|V^{r_1}-v_m|+\big|V^{d}-v^*\big|\big).
\end{eqnarray}
Thus by the wave interaction estimation Lemma 3.1, we have
$$\|G_1\|_{L^1}=O(1)  \d^{ \frac{1}{8}} ( 1+t)^{ -\frac{13}{16}}.$$
Similarly, $\| H_1\|_{L^1}=O(1)  \d^{ \frac{1}{8}} ( 1+t)^{
-\frac{13}{16}}$ can be obtained.

Now we estimate $\| G_2\|_{L^1}$ and $\| H_2\|_{L^1}$. Note that in
$G_2$, besides the wave interaction terms, there are the error terms
due to the $i-$rarefaction waves $(i=1,3).$ So we can write $G_2$ as
$$
\begin{array}{ll}
G_2&\di=- \m
\bigg(\frac{U_{\x}}{V}-\frac{U^b_{\x}}{V^b}-\frac{U^d_{\x}}{V^d}-\sum_{i=1,3}\frac{U^{r_i}_{\x}}{V^{r_i}}\bigg)
_{\x}-\m\bigg(\sum_{i=1,3}\frac{U^{r_i}_{\x}}{V^{r_i}}\bigg)_\x\\
&\di:=G_{21}+G_{22}.
\end{array}
$$
Since the wave interaction terms $G_{21}$ can be verified similarly
as $G_1$, we only compute the error terms $G_{22}$ due to
rarefaction waves.
$$
\begin{array}{ll}
\di \|G_{22}\|_{L^1}&\di =
O(1)\sum_{i=1,3}(\|U^{r_i}_{\x\x}\|_{L^1}+\|(U^{r_i}_\x,V^{r_i}_\x)\|^2)\\
&\di=O(1) \d^{ \frac{1}{8}} ( 1+t)^{ -\frac{13}{16}}
\end{array}
$$
if we choose $q\geq 14$ in Lemma 2.3.

In $H_2$,  besides the wave interaction terms and the error terms
due to the $i-$rarefaction waves $(i=1,3)$, there exists the error
terms $H^d$ due to the viscous $2-$contact wave. We can compute that
$$
\begin{array}{ll}
\di \|H^d\|_{L^1}&\di =
O(1)\d_d(1+t)^{-2}\int_{\mathbf{R}_+}\exp{\left(-\f{C_d(\x+\s_-t)^2}{1+t}\right)}d\x\\
&\di=O(1) \d ( 1+t)^{ -\frac{3}{2}}.
\end{array}
$$
The estimation of $\|G\|$ and $\| H\|$ can be done similarly,  thus
the details are omitted. $\blacksquare$

\subsection{Reformulation of the Problem}

Put the perturbation $(\phi,\psi,  \vartheta)(t,\x)$ around the
superposition wave $(V, U, \T)(t,\x)$  by
\begin{eqnarray}
(\phi,\psi, \vartheta)(t,\x)=(v, u, \theta)(t,\x)-(V, U, \T)(t,\x),
\end{eqnarray}
then by \eqref{(1.7)} and \eqref{(3.1)}, the system for the
perturbation $(\phi,\psi, \vartheta)(t,\x)$ becomes
\begin{eqnarray}\label{31}
\begin{cases}
  \phi_{t}-\s_- \phi_{\x}- \psi_{\x}=   0,  \cr
   \psi_{t}-  \s_- \psi_{\x} + (p-P)_{\x}=
  \m\Big(\frac{u_{\x }}{v}-\frac{U_{\x }}{V}\Big)_\x-G,  \quad\quad  \quad\quad ~~   ~~  t>0,~ \x>0,\cr
   \frac{R}{\g-1}  \big(\vartheta_{t}- \s_-  \vartheta_{\x}\big) + \big(p u_\x-P
   U_\x\big)
   = \k  \left( \frac{\t_{\x }}{v}- \frac{\T_{\x
   }}{V}\right)_\x+\m\Big(
 \frac{(u_{\x})^2}{v}-\frac{(U_{\x})^2}{V}\Big)-H,   \cr (\psi_{0},\psi_{0},\vartheta_0)(\x)
 := (\phi,\psi, \vartheta)(   0, \x)\rightarrow (0,0,0),
  ~~  \text{as} ~ ~ \x \rightarrow+\infty,  \cr (\phi, \psi,\vartheta)(  t,
\x=0)=(V^{d},U^{d},\T^{d})(t,\x=
    0)-(v_m, u_m, \t_m).
\end{cases}
\end{eqnarray}
Define the solution space $\mathbf{X}(0,T)$ to the above system by
\begin{eqnarray}
\mathbf{X} ( 0, T)&:=&\Big\{ ~(\phi,\psi, \vartheta)(t, \x)\,\Big |
\,(\phi,\psi, \vartheta)\in C\left([ 0, T];H^{1} \right)
,~\phi_{\x}\in L^{2}\left( 0, T;L^{2}\right),\cr&&~~~
\big(\psi_{\x}, \vartheta_{\x}\big)\in L^{2}\left(0,
T;H^{1}\right),~ N(T) =:\sup_{0\leq t\leq T} \|(\phi, \psi,
\vartheta)(t)\|_1 \leq \varepsilon_0 \Big\},
\end{eqnarray}
Here $\varepsilon_0\leq \frac{1}{4}\min \bigg\{ \inf\limits_{
\mathbf{R}_+\times\mathbf{R}_+} V(t, \x), \inf\limits_{ \mathbf{R}_+
\times\mathbf{R}_+}
 \T(t, \x)\bigg\}$  is a suitably small and positive constant to be determined.

 Since the proof for the local existence of the solution to the system $(\ref{31})$
 is standard, the details are omitted. To prove
 Theorem 2.1, it is sufficient to prove the  following \emph{a priori} estimate by
 combining the local existence of the solution and the continuation process.

 \vskip 2mm
 \noindent\textbf{Proposition 3.1 (\emph{A priori} estimate)} ~
\emph{Let $(\phi, \psi, \vartheta)\in \mathbf{X} (0,T)$ be a
solution to the system $(\ref{31})$  in the time interval $[0, T)$
with suitably small $\v_0$, and the conditions in Theorem 2.1 hold.
Then there exist a positive constant $C $ independent of $T$  such
that }
\begin{eqnarray}
&&\|(\phi , \psi ,  \vartheta)(t)\|^{2}_{1} +\int^{t}_{0}
\left[\|\phi_{\x}(\tau)\|^{2}+ \|(\psi_{ \x}, \vartheta_\x)
(\tau)\|^{2}_{1} \right] \,d\tau
   \cr&& {}+\int^{t}_{0}\|\sqrt{( U^{b}_{\x},U^{r_1}_{\x},U^{r_3}_{\x})}(\phi,\vartheta)(\tau)\|^2 d\tau \leq  C\left(
\|(\phi_0, \psi_0,  \vartheta_0 )\|^2_1+ \d^{\frac{1}{6}} \right)
   .
\end{eqnarray}

\subsection{Energy estimates}

To prove  Proposition 3.1, we need the following several lemmas.
First we give the following boundary estimates whose proof can be
found in \cite{Qin-Wang}.

\

\noindent{\bf Lemma 3.3 (Boundary Estimates)\cite{Qin-Wang}}~
\emph{There exists the positive constant $C$ such that for any}
$t>0$,
\begin{eqnarray}\nonumber
 &&\int^t_0 |(\phi, \psi, \vartheta)  (\tau, 0) |^2\, d\tau \leq
 C \d, \\\nonumber
 && \int^{t}_{0}\big(\big|\psi\psi_{\x }  \big| +   \big |\vartheta\vartheta_{\x
 } \big|\big)(\tau, 0)\, d\tau   \leq     C \d    +C\d  \int^{t}_{0}
\|( \psi_\x, \vartheta_\x) (\tau ) \|_{1}^2 d\tau.
\\ \nonumber &&
\int^{t}_{0}  \big(|\phi_\tau\psi
  |+  \phi_\x^2\big)(\tau, 0)  \,
 d\tau  \leq C \d + \epsilon\int^{t}_{0} \|\psi_{\x\x}(\tau)\|^2d\tau+  C_\epsilon \int_0^t\|\psi_\x(\tau)\|^2
 d\tau , \\ \nonumber  &&
  \int^{t}_{0} \big( \big|\psi_\tau \psi_\x \big | + \psi_\x^2\big)(\tau, 0)\,
 d\tau  \leq C \d +  \epsilon\int^{t}_{0} \|\psi_{\x\x}(\tau)\|^2d\tau+  C_\epsilon \int_0^t\|\psi_\x(\tau)\|^2
 d\tau, \\ \nonumber  &&
  \int^{t}_{0} \big( \big|\vartheta_\tau \vartheta_\x \big | + \vartheta_\x^2\big)(\tau, 0)\,
 d\tau  \leq C \d+ \epsilon\int^{t}_{0} \|\vartheta_{\x\x}(\tau)\|^2d\tau+  C_\epsilon \int_0^t\|\vartheta_\x(\tau)\|^2
 d\tau,
\end{eqnarray}
\emph{where $\epsilon>0$ is a constant to be determined and
$C_\epsilon$ is the constant depending on $\epsilon$}.

 \vskip 2mm
 \noindent\textbf{Lemma 3.4}~
\emph{Let $(\phi, \psi, \vartheta)\in \mathbf{X} (0,T)$ be a
solution to the system $(\ref{31})$ for some positive T and suitably
small $\v_0>0$, and the conditions in Theorem 2.1 hold. Then there
exist a positive constant $C $ such that}
\begin{eqnarray}\label{(3.22)}
&&\| (\phi, \psi,\vartheta)  (t)\|^{2}_1 + \int^{t}_{0} \| \p_{ \x}
(\tau)\|^{2}+\| ( \psi_\x,\vartheta_\x) (\tau)\|^{2}_1 d\tau \cr
&&{} +\int^{t}_{0}\|\sqrt{(
U^{b}_{\x},U^{r_1}_{\x},U^{r_3}_{\x})}(\phi,\vartheta)(\tau)\|^2
d\tau \cr &\leq & C\left( \|(\phi_0, \psi_0,  \vartheta_0  )\|^2_1+
 \d^{\frac{1}{6}}\right)
   +  C  \d^\frac{1 }{8 }  \int^{t}_{0}
  (1+\tau )^{-\frac{13}{ 12}}
\| (\p, \psi, \vartheta) (\tau)\| ^2 d\tau \cr& &{} +
C\d\int^{t}_{0} \int_{\mathbf{R}_+}
 (1+\tau)^{-1} \exp \left(-\frac{C_d (\x+\s_-  \tau)^2}{1+\tau}
\right)  |(\phi,\vartheta)| ^ 2     d\x d\tau  .
\end{eqnarray}
\noindent \emph{Proof}. \textbf{Step 1.} ~Define
 \begin{eqnarray}
\Phi(\eta):= \eta-\ln \eta-1.
\end{eqnarray}
Under the a priori assumption, there exist  a positive constant $C$
such that
  \begin{eqnarray}
C^{-1}\eta^2\leq\Phi(\eta) \leq C\eta^2.
\end{eqnarray}
Let
\begin{eqnarray}
&&E : =   R \Theta
 \Phi\left(\frac{v}{V} \right) +  \frac{1}{2}\psi^2+
  \frac{R}{\g-1} \Theta
 \Phi\left(\frac{\theta}{\Theta} \right),  \cr&& F:=\s_- E+(P- p)\psi+
\m\bigg(\frac{u_{\x }}{v}-\frac{U_{\x }}{V}\bigg)\psi+
\k\bigg(\frac{\t_{\x }}{v}-\frac{\T_{\x
}}{V}\bigg)\frac{\vartheta}{\t}.
\end{eqnarray}
Then a complicated but direct computation gives
\begin{eqnarray}\label{32}
 E_{t}- F_\x  + \frac{\m\T}{v\theta}\psi_{\x}^2+\frac{\k\Theta}{v
\theta^2}\vartheta_{\x}^2 +P( U^b_\x+U^{r_1}_\x+U^{r_3}_\x) \left[
\g \Phi\left(\frac{v}{V} \right) + \Phi\left(\frac{\t V}{v\T}\right)
\right] =Q ,
\end{eqnarray}
where
\begin{eqnarray}
 Q&=
& -  P U^d_\x \left[  \g \Phi\left(\frac{v}{V} \right) +
\Phi\left(\frac{\t V}{v\T}\right) \right] -\bigg(G\psi  +H  \frac{
\vartheta}{ \t}\bigg)
 \cr& &
+ \bigg[\frac{\m U_{\x } \phi \psi _{\x }}{vV}   + \frac{ 2 \m
U_{\x}\vartheta \psi_{\x}}{v \t } +  \frac{\k\T \T_{\x }\phi
\vartheta_\x}{vV\t^2} + \k \frac{\T_{\x }\vartheta \vartheta_{\x
}}{v\t^2} - \frac{\m (U_{\x})^2 \phi\vartheta}{vV\theta}- \frac{
\k(\T_\x)^2 \phi \vartheta }{vV\t^2}\bigg]\cr& &+ \bigg[\k\left(
\frac{\T_{ \x}}{ V}\right)_\x +\m\frac{(U_{\x})^2}{V}+H
\bigg]\left[(\g-1) \Phi\left(\frac{v}{V} \right) +
\Phi\left(\frac{\theta}{\Theta} \right)-\frac{ \vartheta^2}{
\T\theta} \right] \cr&=:&  \sum_{i=1}^{i=4} Q_i.
\end{eqnarray}
Integrating $(\ref{32})$ over $[0,t]\times \mathbf{R}_+$  yields
\begin{eqnarray}
&&\|(\p, \psi , \vartheta )\|^2  + \int^{t}_{0} \|( \psi_\x,
\vartheta_\x) (\tau ) \|^2
  d\tau+\int_0^t\|\sqrt{(
U^{b}_{\x},U^{r_1}_{\x},U^{r_3}_{\x})}(\phi,\vartheta)(\tau)\|^2
d\tau \cr& \leq&
  C \|(\phi_0, \psi_0,  \vartheta_0 )\|^2+C\int_0^t|F(\tau,\x=0)|d\tau + \sum_{i=1}^{i=4}
I _i
   ,
\end{eqnarray}
where $\di I _i = O(1)\int^{t}_{0} \int_{\mathbf{R}_+} Q _i \,d\x
d\tau$.

\noindent From the boundary estimates in Lemma 3.3, we have
\begin{equation}
\int_0^t|F(\tau,\x=0)|d\tau\leq   C\d
    +C\d \int^t_0
\|( \psi_\x, \vartheta_\x) (\tau ) \|_{1}^2 d\tau.\label{(3.29)}
\end{equation}
We can compute that
\begin{eqnarray}\label{(3.30)}
 I_1\leq
 C \d\int^{t}_{0} \int_{\mathbf{R}_+}
 (1+\tau)^{-1} \exp \left(-\frac{C_d (\x+\s_-\tau)^2}{1+\tau}
\right)   |(\phi,\vartheta)| ^ 2 d\x d\tau
\end{eqnarray}
and
\begin{eqnarray}\label{(3.31)}
   I_2   &\leq & C
   \int^{t}_{0}\|( \psi, \vartheta) (\tau ) \|_{L^\infty}
   (\| G(\tau )\|_{L^1}+\| H(\tau )\|_{L^1})\, d\tau
   \cr&\leq & C \d^{ \frac{1}{8}}
\int^{t}_{0}(1+\tau)^{ -\frac{13}{16} }\|( \psi_\x, \vartheta_\x)
(\tau)\|^{\frac{1}{2}}\|(\psi, \vartheta)(\tau) \|^ {\frac{1}{2}}
d\tau \cr&\leq & \epsilon
   \int^{t}_{0}\|(\psi_\x, \vartheta_\x) (\tau )\|^2 d\tau+C_\epsilon \d^{\frac{1}{6}}
\bigg(1+\int^{t}_{0}( 1+\tau)^{  -\frac{13}{12}}\|( \psi, \vartheta)
(\tau ) \|^2d\tau \bigg)
\end{eqnarray}
where and in the sequel $\epsilon>0$ is a small constant to be
determined and $ C_\epsilon$ is the positive constant depending on
$\epsilon$.

 Now we calculate $I _3$. By Cauchy
inequality, we have
\begin{equation}\label{(3.32)}
I_3\leq
\epsilon\int_0^t\|(\psi_\x,\vartheta_\x)\|^2d\tau+C_\epsilon\int_0^t\int_{\mathbf{R}_+}|(U_\x,\Theta_\x)|^2\cdot|(\phi,\vartheta)|^2d\x
d\tau.
\end{equation}
By Lemma 2.1-Lemma2.3, one has
 \begin{eqnarray}\label{(3.33)}
|(U_\x,\T_\x)|^2
   \leq C  \bigg[\d^{\frac{1}{2}}(1+ t)^{-\frac{3}{2}}+\frac{\d^4}{(1+\d \x)^4}
   + \d(1+t)^{-1} \exp \left(-\frac{C_d (\x+\s_-  t)^2}{1+t}
\right)    \bigg].
   \end{eqnarray}
By the techniques in \cite{Nikkuni-Kawashima}
$$
\begin{array}{ll}
|f(t,\x)|&\di =|f(t,\x=0)+\int_0^\x f_\x(t,\x)d\x|\\
&\di \leq |f(t,\x=0)|+\sqrt{\x}\|f_\x\|,
\end{array}
$$
we have
\begin{eqnarray}\label{(3.34)}
    &&\di \int^{t}_{0}
\int_{\mathbf{R}_+}\frac{\d^4}{(1+\d \x)^4}
   |(\phi,\vartheta)|^ 2 d\x d\tau\cr &\leq&
   C\d^3\int^{t}_{0}|(\phi,\vartheta)(\tau,\x=0)|^2d\tau+C\int_0^t\left[
   \|(\phi_\x,\vartheta_\x)\|^2\int_{\mathbf{R}_+}\frac{\d^4\x}{(1+\d \x)^4}d\x\right] d\tau\cr
   &\leq
   &C\d\left(1+\int_0^t\|(\phi_\x,\vartheta_\x)(\tau)\|^2d\tau\right).
   \end{eqnarray}
Substituting \eqref{(3.33)} and \eqref{(3.34)} into \eqref{(3.32)}
yields
\begin{eqnarray}\label{(3.35)}
   I_3   &\leq&C(\epsilon+ \d)
   \int^{t}_{0}\|(\psi_\x, \vartheta_\x) (\tau ) \| ^2
   d\tau+C\d\int_0^t\|\phi_\x(\tau)\|^2d\tau
 \cr& &{} + C\d
   +C \d ^{\frac{1}{2}}
\int^{t}_{0}(1+\tau)^{-\frac{3}{2}}\|( \phi, \vartheta) (\tau ) \|
^2d\tau  \cr &&{}+ C \d\int^{t}_{0} \int_{\mathbf{R}_+}
 (1+\tau)^{-1} \exp \left(-\frac{C_d (\x+\s_-  \tau)^2}{1+\tau}
\right)   |(\phi,\vartheta)|^2 d\x d\tau
   .
   \end{eqnarray}
Then we have
\begin{eqnarray}\label{(3.36)}
   I_4  &=&O(1)\int_0^t\int_{\mathbf{R}_+} |(\T_{\x\x},V_\x^2,U_\x^2,\T_\x^2,H)||(\phi,\vartheta)|^2d\x
   d\tau.
\end{eqnarray}
So $I_4$ can be estimated similarly as $I_2$ and $I_3$.

Combining \eqref{(3.29)}, \eqref{(3.30)}, \eqref{(3.31)},
\eqref{(3.32)}, \eqref{(3.35)} and \eqref{(3.36)}, and then choosing
$\d$ and $\epsilon$ suitably small yield that
\begin{eqnarray}\label{81}
&&\|(\p, \psi , \vartheta ) (t)  \| ^2  + \int^{t}_{0} \|( \psi_\x,
\vartheta_\x) (\tau ) \| ^2
  d\tau+\int^{t}_{0}\|\sqrt{(
U^{b}_{\x},U^{r_1}_{\x},U^{r_3}_{\x})}(\phi,\vartheta)(\tau)\|^2
d\tau \cr& \leq&
  C   \left(\|(\phi_0, \psi_0,  \vartheta_0 )\|^2 +   \d^\frac{1}{6}
  \right)
    +C\d^{\frac{1}{8}} \int^t_0
\| (  \phi_\x, \psi_{\x\x}, \vartheta_{\x\x}) (\tau ) \| ^2 d\tau
\cr &&{} +C \d ^{\frac{1}{8}}
\int^{t}_{0}(1+\tau)^{-\frac{13}{12}}\|( \phi,\psi, \vartheta) (\tau
) \| ^2d\tau \cr &&{}+ C \d\int^{t}_{0} \int_{\mathbf{R}_+}
 (1+\tau)^{-1} \exp \left(-\frac{C_d (\x+\s_-  \tau)^2}{1+\tau}
\right)   |(\phi,\vartheta)|^ 2   d\x d\tau.
\end{eqnarray}

\noindent\textbf{Step 2.}  ~Differentiating $(\ref{31})_1$ w.r.t.
$\x$ and multiplying it by $\frac{\phi_\x}{v^2 }$  yield
\begin{eqnarray}\label{34}
  \left(\frac{\phi_\x^2}{2v^2}\right)_t
-\s_- \left(\frac{\phi_\x^2}{2v^2}\right)_\x+\frac{ u_x\phi_\x^2 }{
v^3}  - \frac{\phi_\x\psi_{\x\x}}{v^2}=0
  .
\end{eqnarray}
Multiplying $(\ref{31})_2$  by $\frac{\phi_\x}{v}$ gives
\begin{eqnarray}\label{35}
  &&\left(\frac{\phi_\x\psi}{v }\right)_t-\left(\frac{\phi_t\psi}{v
}\right)_\x + \frac{(p-P)_{\x}\phi_\x}{v}\cr &=& - \frac{  U_\x
\phi_\x\psi}{v^2 }+ \frac{V_\x \psi \psi_{\x} }{v^2 } +\s_- \frac{
\phi_{\x}\psi_{\x}}{v}+
 \m\left(\frac{u_{\x}}{v}-\frac{U_{\x}}{V}\right)_\x \frac{\phi_\x}{v}
  -G\frac{
\phi_\x}{v}.
\end{eqnarray}
$\m\times (\ref{34})-(\ref{35})$ gives
\begin{eqnarray}\label{70}
&&\left( \frac{\m\phi_\x^2}{ 2v^2}-\frac{\phi_\x\psi}{v }\right)_t
-\left( \frac{\s_-  \m\phi_\x^2}{2v^2}-\frac{\phi_t\psi}{v }
\right)_\x-\frac{p_v}{v}\phi_\x^2 \cr&=& \frac{  U_\x
\phi_\x\psi}{v^2 }- \frac{V_\x \psi \psi_{\x} }{v^2 } -\s_- \frac{
\phi_{\x}\psi_{\x}}{v}+ \m\frac{ V_\x \phi_\x\psi_\x }{
v^3}-\m\frac{ U_\x \phi_\x^2 }{ v^3}+\m \bigg(\frac{U_\x \phi }{vV}
\bigg)_{\x}\frac{ \phi_\x}{v} \cr& & {}+ \frac{
 p_\t \phi_\x\vartheta_\x}{v}+ \frac{ V_\x(
p_v-P_V)\phi_\x}{v}+ \frac{ \T_\x( p_\t-P_\T)\phi_\x}{v} +G\frac{
\phi_\x}{v}
    .
\end{eqnarray}
Integrating $(\ref{70})$ over $[0,t]\times \mathbf{R}_+$, using the
boundary estimations in Lemma3.3 and choosing $\d$ suitably small
yield
\begin{eqnarray}\label{82}
 &&   \|\phi_\x(t)\|^2
+\int^{t}_{0}  \|\phi_\x(\tau)\|^2
 d\tau \cr&\leq & C\left(
\|(\psi_0, \phi_{0 \x} )\|^2+
 \d ^{\frac{1}{6}}\right)+  C   \d^\frac{1 }{8 }  \int^{t}_{0}
  (1+\tau )^{-\frac{13}{ 12}}
\| (\p, \psi, \vartheta) (\tau)\| ^2 d\tau \cr& & {}
+\int^{t}_{0}\Big\{ C  \Big(\d^{\frac{1}{8}}+\epsilon
  \Big)
 \| ( \psi_{\x\x}, \vartheta_{\x\x})(\tau)\|^2 +
 C_\epsilon
\|\psi_\x(\tau)\|^2\Big\}\,d\tau
 \cr& & {}   +C \d\int^{t}_{0} \int_{\mathbf{R}_+}
 (1+\tau)^{-1} \exp \left(-\frac{C_d (\x+\s_-  \tau)^2}{1+\tau}
\right) |(\phi,\vartheta)|^2 d\x d\tau
     .
\end{eqnarray}

\noindent  \textbf{Step 3.}  ~ Multiplying
 $(\ref{31})_2$
by $-  \psi_{\x\x} $,   then
\begin{eqnarray}\label{36}
 \left(\frac{\psi_\x^2}{2}\right)_t-\left(\psi_{t}\psi_\x-\frac{\s_- }{2} \psi_\x^2
 \right)_\x +  \m\frac{\psi_{\x\x}^2}{v}= \bigg[(p-P)_\x
 + \frac{\m v_{\x}\psi_\x }{v^2}
 +\m \left(\frac{U _{\x }\phi }{vV} \right)_\x+ G \bigg]\psi_{\x\x}.
\end{eqnarray}
Integrating $(\ref{36})$   over $[0,t]\times \mathbf{R}_+$ yields
\begin{eqnarray}\label{83}
&&\|  \psi_\x  (t)\|^{2} + \int^{t}_{0}
 \| \psi_{ \x\x}
(\tau)\|^{2} d\tau  \cr &\leq &C\left( \|(\phi_0, \psi_0,
\vartheta_0  )\|^2_1+
 \d^{\frac{1}{6}}\right)++  C   \d^\frac{1 }{8 }  \int^{t}_{0}
  (1+\tau )^{-\frac{13}{ 12}}
\| (\p, \psi, \vartheta) (\tau)\| ^2 d\tau
  \cr& & {}    +\int^{t}_{0}\Big\{ C  \Big(\d^{\frac{1}{8}}+\epsilon
  \Big)
 \| ( \psi_{\x\x}, \vartheta_{\x\x})(\tau)\|^2 +
 C_\epsilon
\|\psi_\x(\tau)\|^2\Big\}\,d\tau
 \cr& & {}  +C \d\int^{t}_{0} \int_{\mathbf{R}_+}
 (1+\tau)^{-1} \exp \left(-\frac{C_d (\x+\s_-  \tau)^2}{1+\tau}
\right)   \left(\phi^ 2 +\vartheta ^ 2 \right)      d\x d\tau
\end{eqnarray}
where we use the following estimate
\begin{eqnarray}
 \int^{t}_{0} \int_{\mathbf{R}_+} \big|\p_{\x}\psi_\x\psi_{\x\x}\big|\,d\x
d\tau &\leq&
 \int^{t}_{0} \|\phi_\x(\tau)\| \|
  \psi_{\x\x}(\tau)\| \|\psi_\x(\tau)\| _{L^\infty}  d\tau
 \cr&\leq& \int^{t}_{0} \|\phi_\x(\tau)\| \|
  \psi_{\x\x}(\tau)\|^{\frac{3}{2}}  \|\psi_\x(\tau)\| ^{\frac{1}{2}}  d\tau
\cr&\leq& \epsilon \int^{t}_{0} \|
  \psi_{\x\x}(\tau)\|^2   d\tau+C_\epsilon \varepsilon_0^4\int^{t}_{0}
  \|\psi_\x(\tau)\| ^2   d\tau .
 \end{eqnarray}
Multiplying
 $(\ref{31})_3 $
by $-    \vartheta_{\x\x}   $, then
\begin{eqnarray}\label{37}
    &&    \frac{R}{\g-1}  \Bigg[\left(\frac{\vartheta_\x^2}{2}\right)_t-
   \left(\vartheta_{t}\vartheta_\x-\frac{\s_- }{2} \vartheta_\x^2
 \right)_\x\Bigg]+ \frac{\k}{v}\vartheta_{\x\x}^2
 \cr& =
 & \Bigg[\big(pu_\x  - PU_\x\big )
 +\frac{\k v_\x  \vartheta_{\x } }{v^2} +\k \left(\frac{\T _{\x }\phi }{vV} \right)_\x-\m\bigg(
 \frac{(u_{\x})^2}{v}-\frac{(U_{\x})^2}{V}\bigg)+ H \Bigg]\vartheta_{\x\x}
      .
 \end{eqnarray}
Integrating $(\ref{37})$  over $[0,t]\times \mathbf{R}_+$ yields
\begin{eqnarray}\label{84}
&&\| \vartheta_\x  (t)\|^{2} + \int^{t}_{0}
 \| \vartheta_{ \x\x}
(\tau)\|^{2} d\tau  \cr &\leq &C\left( \|(\phi_0, \psi_0,
\vartheta_0  )\|^2_1+
 \d^{\frac{1}{6}}\right)++  C   \d^\frac{1 }{8 }  \int^{t}_{0}
  (1+\tau )^{-\frac{13}{ 12}}
\| (\p, \psi, \vartheta) (\tau)\| ^2 d\tau
  \cr& & {}    +\int^{t}_{0}\Big\{ C  \Big(\d^{\frac{1}{8}}+\epsilon
  \Big)
 \| ( \psi_{\x\x}, \vartheta_{\x\x})(\tau)\|^2 +
 C_\epsilon
\|\vartheta_\x(\tau)\|^2\Big\}\,d\tau
 \cr& & {}   +C \d\int^{t}_{0} \int_{\mathbf{R}_+}
 (1+\tau)^{-1} \exp \left(-\frac{C_d (\x+\s_-  \tau)^2}{1+\tau}
\right)  |(\phi,\vartheta)|^ 2  d\x d\tau,
\end{eqnarray}
where we use the following estimate
\begin{eqnarray}
&&  \int^{t}_{0} \int_{\mathbf{R}_+}
 \big |\p_\x \vartheta_{\x }\vartheta_{\x\x}\big|
 +  \big|\psi_\x^2 \vartheta_{\x\x}\big|\,d\x d\tau \cr&\leq& \epsilon
 \int^{t}_{0} \|(\psi_{\x\x}, \vartheta_{\x\x})(\tau)\|^2
 d\tau+C_\epsilon \varepsilon_0^4 \int^{t}_{0}
   \|(\psi_\x, \vartheta_{\x } )(\tau)\| ^2  d\tau
  .
 \end{eqnarray}
Combining $(\ref{81}), (\ref{82}), (\ref{83})$ and $(\ref{84})$ and
choosing $\d$, $\epsilon$ and $\v_0$ suitably small, we can complete
the proof of Lemma 3.4. $\blacksquare$

 \vskip 2mm

Now to close the a priori estimates, the remaining thing is to
compute the last term
 in the right-hand side of $(\ref {(3.22)})$ which comes from the viscous contact
 wave. Here we use the method of the heat kernel estimation invented
 in \cite{Huang-Li-Matsumura}.

\noindent {\bf Lemma 3.5.\cite{Huang-Li-Matsumura} }  \emph{Suppose
that $h(t,\x)$ satisfies
\begin{eqnarray}
h\in L^\infty\left(0,T; L^2(\mathbf{R}_+)\right),~~h_\x\in
L^2\left(0,T; L^2(\mathbf{R}_+)\right),~~h_t-\s_-  h_\x\in
L^2\left(0,T; H^{-1}(\mathbf{R}_+)\right),
\end{eqnarray}
then
\begin{eqnarray}\label{(3.50)}
 &&\int_0^t
\int_{\mathbf{R}_+}(1+\tau)^{-1}\exp
\left(-\frac{2a(\x+\s_-\tau)^2}{1+\tau} \right) h^2 \,d\x d\tau
\cr&\leq & C_a\left\{ \| h(0,\cdot ) \|^2+\int_0^t\Big [ h^2(\tau,
0) +\| h_\x (\tau,\cdot)\|^2 +\big{\langle} h_\tau-\s_- h_\x,  (w^a
)^2 h \big{\rangle}_{H^{-1}\times H^{ 1}}\Big ] d\tau\right\}
 \end{eqnarray}
where
 \begin{eqnarray}
w^a(t,\x)=-(1+t)^{-\frac {1}{2}}\int_{\x+\s_- t}^{\infty} \exp
\left(-\frac{a y^2}{1+t}   \right)d y,
\end{eqnarray}
and $a>0$ is a constant to be determined}.

\vskip 2mm

Based on Lemma 3.5, we have the desired estimates in the following
Lemma.

\

\noindent{\bf Lemma 3.6}  \emph{There exist a uniform constant $C>0$
such that if $\d$ and  $\varepsilon_0$ are small enough, then we
have
\begin{eqnarray}
 &&\int_0^t
\int_{\mathbf{R}_+}(1+\tau)^{-1}\exp
\left(-\frac{C_d(\x+\s_-\tau)^2}{1+\tau} \right) |(\p,\psi,
\vartheta)|^2 \,d\x d\tau \cr&\leq & C\left( \|(\phi_0, \psi_0,
\vartheta_0  )\|^2_1+
 \d ^{\frac{1}{6}} \right)
        +    C    \d^\frac{1 }{ 8 }\int^{t}_{0}
  ( 1+\tau)^{- \frac{13}{12}}
\|(\phi , \psi , \vartheta)(\tau)\|^{2} d\tau
   .
 \end{eqnarray}}
\noindent \emph{Proof}. \textbf{Step 1.} First, let
\begin{eqnarray}
h = P\phi+ \frac{R}{\g-1} \vartheta
\end{eqnarray}
in Lemma 3.4. Then we only need to control the last term of
\eqref{(3.50)} on the right hand side.

We have from the energy equation $\eqref{31}_3$,
\begin{eqnarray}
h_t-\s_-  h_\x&=&   (P- p)  \psi_\x + U_\x (P-p)  + \big(P_t-\s_-
P_\x\big)\phi     \cr & & {}+  \k  \big( \frac{\t _{ \x}
}{v}-\frac{\T_{ \x}}{ V} \big)_\x
 +  \m \big( \frac{(u_{\x})^2}{v} -\frac{ ( U_\x)^2}{ V}
   \big)-H .
\end{eqnarray}
Thus
\begin{eqnarray}\label{47}
&& \int_0^t\big{\langle} h_\tau-\s_- h_\x, (w^a)^2 h
\big{\rangle}_{H^{-1}\times H^{ 1}}\, d\tau  \cr &=& - \k
\int^{t}_{0}
  \big[(
\frac{\t_{\x }}{v}-\frac{\T_{ \x}}{ V} ) (w^a)^2h \big] (\tau, 0)\,
d\tau -
   \k \int^{t}_{0} \int_{\mathbf{R}_+}  (
\frac{\t_{\x }}{v}-\frac{\T_{ \x}}{ V} )[ (w^a)^2 h]_\x\, d\x d\tau
\cr&  & {}+\int^{t}_{0} \int_{\mathbf{R}_+}\big[(P- p) \psi_\x +
U_\x (P-p) + \big(P_t-\s_- P_\x\big)\phi\cr&&  \qquad\qquad\qquad+\m
\big( \frac{(u_{\x})^2}{v} -\frac{ ( U_\x)^2}{ V}
   \big)  -H \big]\left
(w^a\right)^2h \,d\x d\tau .
\end{eqnarray}
Notice that
\begin{eqnarray}
 && \| w^a(t) \|_{L^\infty}\leq C_a, \quad  w^a _\x =    (1+t)^{-\frac{1}{2}}
 \exp \left(-\frac{ a (\x+\s_-  t ) ^2}{ 1+t }\right),
 \quad \big|w^a _t-\s_- w^a _\x \big|\leq C_a (1+t)^{-1},\cr&&
\big |P_t-\s_-  P_\x\big| \leq     C \bigg\{  U^b_\x +U^{r_1}_\x +
U^{r_3}_\x +
     \d(1+t)^{-1}\exp \left(-\frac{ C_d (\x+\s_-  t ) ^2}{ 1+t
     }\right)
      \bigg\},
\end{eqnarray}
thus to control terms on the right hand side of $(\ref{47})$, we
only consider the term $(w^a )^2(P- p)h \psi_\x$. By using the mass
equation $\eqref{31}_1$  and the momentum equation $\eqref{31}_2$
again, we have
\begin{eqnarray}\label{(3.57)}
  \left (w^a\right)^2(P- p) h \psi_\x
 &=  &      \frac{  (w^a)^2 [\g P \phi-(\g-1)h]h (\p_t-\s_-
 \p_\x)}{v}
 \cr&=&  \frac{ \g  P (w^a)^2 h}{2v} \Big[\big(\p^2\big)_t-\s_- \big( \p^2\big)_\x\Big]
 -\frac{ (\g-1) (w^a)^2  h^2}{v}  \big(\p_t-\s_-  \p_\x\big)
 \cr&=& \left(  \frac{\g P (w^a)^2  h
\p^2 -2(\g-1)(w^a)^2 \p h^2 }{2v}\right)_t\cr& &{}- \s_- \left(
\frac{\g P (w^a)^2  h \p^2 -2(\g-1)(w^a)^2 \p h^2 }{2v}\right)_\x
\cr& &{}-  \frac{\g Ph \p^2 -2(\g-1) \p h^2 }{ v}w^a\big(
 w^a  _t-\s_-   w^a  _\x\big)-\frac{\g (w^a)^2\p^2
h}{2v}\big(P_t-\s_- P_\x\big) \cr& & {}+  \frac{\g P (w^a)^2 h\p^2
-2(\g-1)(w^a)^2 \p h^2 }{2 v^2} \big(  \psi_\x +U_\x\big) \cr& &{} +
\frac{(w^a)^2[ 4(\g-1) h- \g P \p ]\p
  }{2 v } \big ( h_t-\s_-  h_\x\big).
\end{eqnarray}
Now the terms in the right hand side of \eqref{(3.57)} can be
estimated directly and in particular, we have
\begin{eqnarray}
&& \int^{t}_{0} \int_{\mathbf{R}_+} \frac{\g P (w^a)^2 h\p^2
-2(\g-1)(w^a)^2 \p h^2 }{2 v^2}  \psi_\x d\x d\tau \cr&\leq&C
\int^{t}_{0} \int_{\mathbf{R}_+} |\psi_\x| \big ( |\p|^3
+|\vartheta|^3\big) \,d\x d\tau\cr &\leq
&C\int_0^t\|(\phi,\vartheta)\|^2_{L_\i}\|\psi_\x\|\|(\phi,\vartheta)\|d\tau\cr
 & \leq&
C\varepsilon_0^2 \int^{t}_{0} \|  (\p_\x, \psi_\x , \vartheta_\x)
(\tau ) \| ^2 d\tau.
\end{eqnarray}
The other terms can be controlled by the similar procedure as Step 1
of Lemma 3.4. Thus the combination of the above estimates and Lemma
3.5 yield
\begin{eqnarray}\label{48}
 &&\int_0^t
\int_{\mathbf{R}_+}(1+\tau)^{-1}\exp
\left(-\frac{2a(\x+\s_-\tau)^2}{1+\tau} \right) \bigg(P\phi+
\frac{R}{\g-1} \vartheta\bigg)^2 \,d\x d\tau \cr&\leq & C_a\left(
\|(\phi_0, \psi_0,  \vartheta_0  )\|^2_1+
 \d  ^{\frac{1}{6}} \right)
 +    C_a   \d ^\frac{1 }{ 8} \int^{t}_{0}   ( 1+\tau)^{-
\frac{13}{12}} \|(\phi , \psi , \vartheta)(\tau)\|^{2} d\tau
 \cr \cr& & {} +C_a (\d+\varepsilon_0)\int^{t}_{0} \int_{\mathbf{R}_+}
 (1+\tau)^{-1} \exp \left(-\frac{C_d (\x+\s_-  \tau)^2}{1+\tau}
\right)   |(\phi,\vartheta)| ^ 2  d\x d\tau .
 \end{eqnarray}

\noindent\textbf{Step 2.} ~ Let
\begin{eqnarray}
  W^A(t, \x) & := &  -(1+t)^{-1} \int^\infty_{\x+\s_-  t  }
 \exp \left(-\frac{ A y ^2}{ 1+t }\right)
  \,d y ,
\end{eqnarray}
 where $A>0$ is a   constant to be determined.

  Then
\begin{eqnarray}
W^A_\x= (1+t)^{-1}
 \exp \left(-\frac{ A (\x+\s_-  t ) ^2}{ 1+t }\right),
 \quad  \big|W^A_t-\s_- W^A_\x\big|  \leq C_A (1+t)^{-\frac{3}{2}}.
\end{eqnarray}
From the fact $p-P= \frac{R\vartheta-P \phi}{v} $, we have
\begin{eqnarray}\label{+}
  \frac{( R\vartheta-P \phi )_\x}{v}-\frac{v_\x(R\vartheta-P
  \phi)}{v^2 }=- \big(\psi_t-\s_-  \psi_\x\big)+
 \m  \left(\frac{u_{\x }}{v}-\frac{U_{\x }}{V}\right)_\x- G.
\end{eqnarray}
Multiplying \eqref{+} by $W^A( R\vartheta-P \phi )$ implies
\begin{eqnarray}
 && \left(\frac{W^A ( R\vartheta-P \phi )^2}{2v}\right)_\x
 -\frac{W^A_\x ( R\vartheta-P \phi )^2}{2v}
-\frac{ W^A v_\x (R\vartheta-P
  \phi)^2}{2v^2} \cr & = & - W^A \bigg[ \big(\psi_t-\s_-
  \psi_\x\big)-
 \m   \left(\frac{u_{\x }}{v}-\frac{U_{\x }}{V}\right)_\x+ G\bigg]( R\vartheta-P \phi ).
\end{eqnarray}
Note that
$$
\begin{array}{ll}
\di W^A(\psi_t-\s_-  \psi_\x)( R\vartheta-P \phi )
 &\di= \big\{ W^A \psi  ( R\vartheta-P \phi )\big\}_t-\s_- \big \{ W^A \psi  ( R\vartheta-P \phi )\big\}_\x
 \\
&\di  -\psi  ( R\vartheta-P \phi )\big(  W^A _t-\s_-  W^A _\x\big)
\\
&\di- W^A \psi \big\{ (R\vartheta-P \phi )_t-\s_- (R\vartheta-P \phi
) _ \x\big\},
 \end{array}
 $$
\begin{eqnarray}
  &&  (R\vartheta-P \phi
)_t-\s_-  (R\vartheta-P \phi ) _ \x
  \cr&  = & {}  (\g-1)\bigg\{( P-p)  u_\x
   + \k  \left( \frac{\theta_{\x }}{v}-\frac{\T_{\x }}{V}\right)_\x+\m
\bigg( \frac{(u_{\x})^2}{v}- \frac{(U_{\x})^2}{V}\bigg)-H\bigg\}
\cr& &{}-
  \g P \psi_\x-\big(P_t-\s_-  P_\x\big)\phi
 \end{eqnarray}
 and
 \begin{eqnarray}
 \g   P W^A \psi  \psi_\x& =&  \frac{\g}{2} \big(  P W^A \psi^2 \big)
 _\x- \frac{\g}{2}    P W^A _\x\psi^2
 - \frac{\g}{2} P   _\x W^A      \psi^2,
 \end{eqnarray}
we have
\begin{eqnarray}\label{49}
   -\frac{W^A_\x  }{2v}\big\{ ( R\vartheta-P \phi )^2+ \g v   P
   \psi^2\}
    =   -\big\{ W^A \psi  ( R\vartheta-P \phi )\big\}_t- E^{ A}_\x
 +Q^{ A},
 \end{eqnarray}
where
\begin{eqnarray}
 E^{ A}:&=&\frac{W^A ( R\vartheta-P \phi )^2}{2v}+\frac{\g}{2}
  P W^A \psi^2-
 \m   W^A ( R\vartheta-P \phi )\left(\frac{u_{\x } }{v}
 -\frac{ U_\x}{V}\right)   \cr&   &{}-\s_-    W^A \psi  ( R\vartheta-P \phi
 )-(\g-1)\k W^A \psi\left(
\frac{\theta_{\x }}{v}-\frac{\T_{ \x}}{ V} \right)   ,
\end{eqnarray}
and
\begin{eqnarray}
 Q^{ A}:&=&   \frac{ W^Av_\x (P-p )^2}{2 } +\big( W^A _t-\s_-
 W^A _\x\big)   ( R\vartheta-P \phi ) \psi
 - W^AG( R\vartheta-P \phi )
  \cr&   & {}+
W^A\psi\Bigg\{(\g-1) \Bigg[ (P-p) u_\x+\m\Bigg( \frac{u_{\x}^2}{v} -
\frac{ ( U_\x)^2}{ V}
    \Bigg) -H\Bigg]-\big (P_t-\s_-  P_\x\big) \p+\frac{\g P_\x\psi}{2}\Bigg\}
    \cr &   &{}  -\m  \big\{ W^A ( R\vartheta-P \phi )
 \big\}_\x\left(\frac{u_{\x }}{v}
 -\frac{ U_\x}{V}\right)
  -(\g-1)\k \big(W^A \psi \big)_\x\left(
\frac{\t_{\x } }{v}-\frac{\T_{ \x}}{ V}  \right) .
\end{eqnarray}
First, we have
\begin{eqnarray}
 \bigg |\int^t_0  E^{ A}(\tau, 0)
\,d\tau \bigg | \leq   C_A \d    +C_A\d   \int^t_0 \|( \psi_\x,
\vartheta_\x) (\tau ) \|_{1}^2 d\tau .
\end{eqnarray}
The estimations of the terms concerned with  $ W^A$ are similar to
those in Step 1 while the other terms are similar to those of  Step
1 in the proof of Lemma 3.4. Thus integrating  $(\ref{49})$ over
$[0,t]\times \mathbf{R}_+$ yields
\begin{eqnarray}\label{50}
  && \int^{t}_{0}
\int_{\mathbf{R}_+} (1+\tau )^{-1}\exp \left(-\frac{ A (\x+ \s_-
\tau)^2}{ 1+\tau }\right)
   \big\{( R\vartheta-P \phi )^2+   \psi^2 \big\}\, d\x d\tau
 \cr&\leq & C_A\left(
\|(\phi_0, \psi_0,  \vartheta_0  )\|^2_1+
 \d  ^{\frac{1}{6}} \right)
 +    C_A   \d ^\frac{1 }{ 8} \int^{t}_{0}   ( 1+\tau)^{-
\frac{13}{12}} \|(\phi , \psi , \vartheta)(\tau)\|^{2} d\tau
 \cr \cr& & {} +C_A (\d+\varepsilon_0)\int^{t}_{0} \int_{\mathbf{R}_+}
 (1+\tau)^{-1} \exp \left(-\frac{C_d (\x+\s_-  \tau)^2}{1+\tau}
\right)   |(\phi,\vartheta)| ^ 2  d\x d\tau .
\end{eqnarray}

\noindent\textbf{Step 3.} Combining $(\ref{48})$ and $(\ref{50})$,
 then choosing $A=2a=C_d$ and setting
 $\d,\, \varepsilon_0 $ suitably small, we can complete the proof of Lemma 3.6. $\blacksquare$

 \vskip 2mm
 \noindent\emph{Proof of Proposition 3.1.}  Choosing $\d,\v_0$ suitably small in Lemmas 3.4 and
 Lemma 3.6, then using Gronwall inequality yield Proposition 3.1. $\blacksquare$


\begin{thebibliography}{000}
\bibitem{Duan-Liu-Zhao} R. Duan, H. Liu,  H.  Zhao,  Global stability of rarefaction
waves for the compressible Navier-Stokes equations.  Trans. Amer.,
Math. Soc., 361 (2009), no. 1, pp. 453--493.

\bibitem{Huang-Li-Matsumura} F. Huang,  J.  Li, A.  Matsumura,
 Stability of the
combination of the viscous contact wave and the rarefaction wave to
the compressible Navier-Stokes equations. Arch. Rat. Mech. Anal.
(DOI) 10.1007/s00205-009-0267-0.


\bibitem{Huang-Li-Shi} F. Huang,  J.  Li, X.  Shi,  Asymptotic behavior of solutions
to the full compressible Navier-Stokes equations in the half space.
 to appear in Commu. Math. Sci.



\bibitem{Huang-Matsumura-Shi} F. Huang,  A.  Matsumura, X.  Shi, Viscous shock wave
and boundary layer solution to an inflow problem for compressible
viscous gas. Comm. Math. Phys., 239(2003),  pp. 261--285.



\bibitem{Huang-Matsumura-Xin} F. Huang, A.  Matsumura,
Z. Xin,  Stability of contact discontinuities for the 1-D
compressible Navier-Stokes equations. Arch. Rat. Mech. Anal.,
179(2005),
 pp.  55--77.



\bibitem{Huang-Qin} F. Huang,  X. Qin,  Stability of  boundary layer
and rarefaction wave to an outflow problem for  compressible
Navier-Stokes equations  under large perturbation.    J. Diff.
Eqns., 246(2009), pp. 4077--4096.

\bibitem{Huang-Xin-Yang}F. Huang, Z. Xin, T. Yang, Contact
discontinuity with general perturbations for gas motions. Adv.
Math., 219 (2008), no. 4, 1246--1297.


\bibitem{Kawashima-Matsumura} S. Kawashima, A. Matsumura, Asymptotic
stability of traveling wave solutions of systems for one-dimensional
gas motion. Comm. Math. Phys., 101 (1985), 97--127.

\bibitem{Kawashima-Matstumura-Nishihara} S. Kawashima, A. Matsumura, K. Nishihara,
Asymptotic behavior of solutions for the equations of a viscous
heat-conductive gas, Proc. Japan Acad., 62A (1986), 249--252.

\bibitem{Kawashima-Nishibata-Zhu}  S. Kawashima, S. Nishibata, P. Zhu,
Asymptotic stability of the stationary solution to the compressible
 Navier-Stokes equations in the half space. Comm. Math. Phys., 240 (2003), no. 3, 483--500.

\bibitem{Liu} T. Liu, Shock waves for compressible Navier-Stokes
equations are stable. Comm. Pure Appl. Math., XXXIX (1986),
565--594.


\bibitem{Liu-Matsumura-Nishihara} T. Liu, A. Matsumura and K.
Nishihara, Behaviors of solutions for the Burgers equation with
boundary corresponding to rarefaction waves. SIAM J. Math. Anal.,
29(1998), 293-308.



\bibitem{Liu-Xin-1} T. Liu, Z. Xin, Nonlinear stability of rarefaction
waves for compressible Navier-Stokes equations. Comm. Math. Phys.,
118 (1988), 451--465.

\bibitem{Liu-Xin-2} T. Liu, Z. Xin, Pointwise decay to contact
discontinuities for systems of viscous conservation laws.  Asian J.
Math., 1 (1997),  no. 1, 34--84.

\bibitem{Matsumura} A. Matsumura,  Inflow and outflow problems in the half space
 for a one-dimensional isentropic model system of compressible viscous
 gas. Proceedings of IMS Conference
on Differential Equations from Mechanics (Hong Kong, 1999).


\bibitem{Matsumura-Nishihara-1} A. Matsumura, K. Nishihara, On the stability
of traveling wave solutions of a one-dimensional model system for
compressible viscous gas. Japan J. Appl. Math., 2 (1985), 17--25.

\bibitem{Matsumura-Nishihara-2} A. Matsumura, K. Nishihara, Asymptotics toward
the rarefaction wave of the solutions of a one-dimensional model
system for compressible viscous gas. Japan J. Appl. Math., 3 (1986),
1--13.

\bibitem{Matsumura-Nishihara-3} A.Matsumura, K. Nishihara, Large-time behavior
 of solutions to an inflow problem in the half space
for a one-dimensional system of compressible viscous gas. Comm.
Math. Phys., 222 (2001), 449--474.


\bibitem{Nikkuni-Kawashima} Y. Nikkuni, S. Kawashima, Stability of
stationary solutions to the half-space problem for the discrete
Boltzmann equation with multiple collisions. Kyushu J. Math., 54
(2000), no. 2, 233--255.

\bibitem{Nishihara-Yang-Zhao} K. Nishihara, T. Yang,   H. Zhao, Nonlinear stability of
strong rarefaction waves for compressible Navier-Stokes equations.
SIAM J. Math. Anal.,  35  (2004),  no. 6, 1561--1597.

\bibitem{Qin-Wang}  X. Qin, Y. Wang, Stability of wave patterns to the inflow problem of full compressible Navier-Stokes equations.
 SIAM J. Math. Anal., 41 (2009), 2057-2087.

\bibitem{Smoller} J. Smoller, Shock waves and reaction-diffusion equations, Berlin,
Heidelberg, New~York, Springer 1982.

\bibitem{Szepessy-Xin}A. Szepessy, Z. Xin, Nonlinear stability of viscous
shock waves. Arch. Rat. Mech. Anal., 122 (1993), no. 1, 53--103.

\bibitem{Xin} Z. Xin, On nonlinear stability of contact
discontinuities. Hyperbolic problems: theory, numerics, applications
(Stony Brook, NY, 1994), 249--257, World Sci. Publ., River Edge, NJ,
1996.

\bibitem{Zhu} P.  Zhu, Existence and asymptotic stability of stationary solution to the full
compressible Navier-Stokes equations in the half space.
S¨±rikaisekikenky¨±sho K¨­ky¨±roku, No. 1247 (2002), 187--207.

 \end{thebibliography}
\end{document}